\documentclass[ejs]{imsart}
\usepackage{amsfonts}
\usepackage{amsmath}
\usepackage{amssymb,mathrsfs}
\usepackage{amsthm,nicefrac,bigints}
\usepackage[utf8]{inputenc}
\usepackage{epsfig}
\usepackage{multirow}
\usepackage{eurosym}
\usepackage{enumerate}

\newtheorem{theorem}{Theorem} 
\newtheorem{proposition}[theorem]{Proposition}

\theoremstyle{definition}\newtheorem{remark}[theorem]{Remark}
\theoremstyle{definition}\newtheorem{example}[theorem]{Example}

\newcommand{\bB}{\mathbb{B}}

\newcommand{\bM}{\mathbb{M}}
\newcommand{\bN}{\mathbb{N}}

\newcommand{\bQ}{\mathbb{Q}}
\newcommand{\Qbin}{\mathbb{Q}_\text{\rm bin}}

\newcommand{\bR}{\mathbb{R}}

\newcommand{\bZ}{\mathbb{Z}}

\newcommand{\pp}{v}

\newcommand{\cA}{\mathcal{A}}
\newcommand{\cB}{\mathcal{B}}

\newcommand{\cD}{\mathcal{D}}
\newcommand{\cE}{\mathcal{E}}
\newcommand{\cF}{\mathcal{F}}

\newcommand{\cM}{\mathcal{M}}

\newcommand{\cP}{\mathcal{P}}

\newcommand{\cS}{\mathcal{S}}

\newcommand{\fE}{\mathfrak{E}}
\newcommand{\fF}{\mathfrak{F}}
\newcommand{\fM}{\mathfrak{M}}
\newcommand{\fR}{\mathfrak{R}}
\newcommand{\fl}{\mathfrak{l}}

\newcommand{\unif}{\text{\rm unif}}
\newcommand{\var}{\text{\rm var}}

\newcommand{\Po}{\text{\rm Po}}
\newcommand{\Mix}{\text{\rm Mix}}

\newcommand{\Dexp}{\text{\rm DExp}}
\newcommand{\MSE}{\text{\rm MSE}}

\newcommand{\Leb}{\fl} 

\allowdisplaybreaks

\newcommand{\abs}[1]{\lvert#1 \rvert}

\begin{document}

\begin{frontmatter}

\title{Discrete mixture representations of 
	parametric distribution families: geometry and statistics}
\runtitle{Discrete mixture representations}

\begin{aug}
\author{\fnms{Ludwig} \snm{Baringhaus}\ead[label=e1]{lbaring@stochastik.uni-hannover.de}}
\and
\author{\fnms{Rudolf} \snm{Gr{\"u}bel}\ead[label=e2]{rgrubel@stochastik.uni-hannover.de}}
\address[]{Institut f\"ur Mathematische Stochastik, Leibniz Universit\"at Hannover, Welfengarten 1,
 D-30167 Hannover, Germany \printead{e1,e2}}

\runauthor{L. Baringhaus and R. Gr{\"u}bel}

\affiliation{Leibniz Universit\"at Hannover}

\end{aug}

\begin{abstract}
	We investigate existence and properties of discrete mixture 
	representations $P_\theta =\sum_{i\in E} w_\theta(i) \, Q_i$ for a given
	family $P_\theta$, $\theta\in\Theta$, of probability measures. The 
	noncentral chi-squared distributions provide a classical example.
	We obtain existence results and results about geometric and statistical 
	aspects of the problem, the latter including loss of Fisher information,
	Rao-Blackwellization, asymptotic efficiency and nonparametric 
	maximum likelihood estimation of the mixing probabilities.
\end{abstract}



\end{frontmatter}

\section{Introduction}\label{sec:intro}

We say that a family $\{P_\theta:\, \theta\in\Theta\}$ of probability measures on some
measurable space $(S,\cS)$ has a mixture representation in terms of a finite or
countably infinite family  
$\{Q_i:\, i\in E\}$ of probability measures on $(S,\cS)$, the \emph{mixing distributions},
and a family of $\{w_\theta:\, \theta\in\Theta\}$ of probability mass functions on $E$, the
\emph{mixing coefficients} if, for all $\theta\in\Theta$,
\begin{equation}\label{eq:mixdef}
	P_\theta(A) =\sum_{i\in E} w_\theta(i) \, Q_i(A)\quad\text{for all } A\in\cS.
\end{equation}
We will generally have a continuous (uncountable) base family $\{P_\theta:\, \theta\in\Theta\}$
and a parameter set $\Theta$ that is a subset of $\bR^d$ for some $d\ge 1$. 
Note that we assume that $E$ is discrete; we will therefore refer 
to~\eqref{eq:mixdef} as a \emph{discrete mixture representation}.

A particularly interesting example, with connections to the power of statistical tests and
also to the Dynkin isomorphism in the theory of stochastic processes, arises if we start with a
standard normal random variable $X$ and let $P_\theta$ be the distribution of $(X+\theta)^2$.
For such noncentral chi-squared distribution with one degree of freedom and
noncentrality parameter~$\theta^2$, the representation~\eqref{eq:mixdef} holds with $E=\bN_0$,
$Q_i$ the central chi-squared distribution with $2i+1$ degrees of freedom and $w_\theta$ the 
probability mass function of the Poisson distribution with mean $\theta^2/2$; see 
also equation~\eqref{eq:chinrepr} below. 

Replacing the standard normal distribution by some other 
distribution $\mu$ on the real line 
we obtain the noncentral distributions associated with $\mu$. In~\cite{BaGr5} we
obtained representations similar to the classical case for such noncentral families if
$\mu$ is the logistic distribution or the hyperbolic secant distribution, 
but our methods there did not lead to similar results for
other standard symmetric distributions, such as the double exponential and the Cauchy distribution.
Initially, our aim was to develop an alternative approach that can be used in these cases and
for other parametric families.  It turned out that
the problems of existence and properties of discrete mixture representations for a given family of
probability measures have some interesting general aspects, of a geometric and statistical 
nature; these, together with their interaction, are now the major themes of the present paper.

We briefly recall the general situation: With an arbitrary set of mixing distributions
we need to replace the sum in~\eqref{eq:mixdef} by an integral and thus require a measurable structure,
i.e.\ a $\sigma$-field~$\cE$, on $E$ together with $\cE$-measurability of all 
functions $y\mapsto Q(y,A)$, $A\in \cS$, so that $Q$ is a Markov kernel (transition  probability)
from $(E,\cE)$ to $(S,\cS)$.  For a probability
measure $\tau$ on $(E,\cE)$ we then obtain $P_\tau$, the mixed distribution
with mixing kernel $Q$ and mixing measure $\tau$,  by
\begin{equation}\label{eq:mixgen}
   P_\tau(A) := \int Q(y,A)\, \tau(dy)\quad\text{for all } A\in \cS.
\end{equation} 
With $S=\{y_k:\, k\in \bN_0\}$ countable and $\cS$ the set of all subsets of $S$ a discrete 
mixture representation always exists for any family $\{P_\theta:\, \theta\in\Theta\}$: 
With $E=\bN_0$ we take $Q_k$ to be the probability measure $\delta_{y_k}$ concentrated 
on $\{y_k\}$ and define 
the mixing coefficients by $w_\theta(k)=P_\theta(\{y_k\})$. In the general situation we may 
still use such a construction with 
$Q(y,\cdot)=\delta_y$ the one-point measure in $y$ and $P$ itself as the mixing measure
(here $E=S$, $\cE=\cS$,  and we assume that $\{y\}\in \cS$ for all $y\in S$).
In particular there is a mixture representation for \emph{every} 
family $\{P_\theta: \, \theta\in\Theta\}$
of probability measures on a given measurable space $(S,\cS)$. A \emph{discrete} mixture 
representation, however, might not exist: For example, \eqref{eq:mixdef} implies that the family 
$\{P_\theta:\, \theta\in\Theta\}$ is dominated by some $\sigma$-finite measure $\nu$
(this is not true in general, for example if $E$ is uncountable and $E=S$, and if the family 
consists
of all one-point measures). 
In fact, for a dominated family~\eqref{eq:mixdef} may be rewritten in terms of
functions as
\begin{equation}\label{eq:mixdef0}
   f_\theta(x) =\sum_{i\in E} w_\theta(i) \, g_i(x)
                        \quad\text{for $\nu$-almost all } x\in S,
\end{equation}
where $g_i,f_\theta$ are $\nu$-densities of $Q_i,P_\theta$. This way, we may regard
$\{P_\theta:\, \theta\in\Theta\}$ as a curve in the space $L_1=L_1(S,\cS,\nu)$ of 
$\nu$-integrable functions on $S$. 

Mixed distributions are a canonical theme in probability and statistics, and
many authors have considered related problems.
A standard reference for mixtures, especially from a statistical point of view, is 
Lindsay's research monograph~\cite{Lindsay}. 
The set of mixtures of binomial distributions has been considered in detail 
in~\cite{WoodProb} and~\cite{WoodStat}, from a geometric and statistical
angle respectively. In~\cite{Leroux} penalized maximum-likelihood 
estimators are proposed for mixed distributions.
A classical mixture result for a nonparametric distribution family
appears in connection with Grenander's influential paper~\cite{Grenander} 
about estimation of  distributions on $\bR_+$ that have decreasing densities
and, more generally, in connection with the structure of unimodal 
distributions~\cite[p.158]{FellerII}.  
The special case of noncentral chi-squared is considered in many papers, see
e.g.~\cite{GoelDeGroot79,SaxenaAlam82,Kubokawa17}.
Finally, a strong case for the use of mixture representations is made in~\cite{Hoff}, 
where these are related to the removal of constraints in the optimization problems 
that typically turn up once the estimates have to be calculated.

In Section~\ref{sec:existence} we provide existence results that can be used to obtain discrete
mixture representations in many cases, including the noncentral families associated with the
double exponential and the Cauchy distribution mentioned above; 
see Theorems~\ref{thm:existunif} and \ref{thm:noncentral1}. The proof of the first theorem
is based on a representation of the $\sigma$-field $\cS$ by a filtration that consists of
$\sigma$-fields generated by finite or countably infinite partitions of $S$, 
such as the dyadic partitions if $\cS$ is 
the $\sigma$-field of Borel subsets of the unit interval $S$. 
The result, however, is less explicit than in the
classical case of noncentral chi-squared distributions or the other cases considered in~\cite{BaGr5}. 

In Section~\ref{sec:continuity} we develop the curve
view in~\eqref{eq:mixdef0} and relate the existence of discrete mixture representations
to continuity properties of the curve for spaces other than~$L_1$. 

Section~\ref{sec:geo} deals with geometric aspects.
The set of \emph{all} mixtures of a given family of mixing distributions is obviously convex,
but even in a stronger sense which makes it well-suited to Dynkin's approach~\cite{Dynkin78} 
that is based 
on the notion of barycentric map. In particular, and in contrast to the strategy used by
many papers in this area where Choquet's representation theorem is an important tool, we do not 
require topological notions for infinite-dimensional linear spaces. 
Naturally, the (geometric) question
arises whether a representation such as~\eqref{eq:mixdef} is minimal in the sense
that the right hand side is the barycentric convex hull of the parametric family on the
left; see Theorem~\ref{thm:simplexchi} for results in the chi-squared case. 
The set of extreme points of the mixture family in the representation of uniforms from
Section~\ref{sec:existence} turns out to be empty; see part (a) of Theorem~\ref{thm:unifconv}.

In Section~\ref{sec:stat}, which contains several subsections,  
we consider statistical aspects. In particular, from 
this point of view mixtures such as~\eqref{eq:mixdef} may be thought of as 
representing a two-stage experiment: To obtain a random variable $X$ with distribution
$P_\theta$, we first choose an $E$-valued random variable $T$ with probability mass function
$w_\theta$ and then, given $T=i$, choose $X$ with distribution $Q_i$. In particular, finding
a discrete mixture representation is essentially the same as finding a discrete sufficient statistic,
on a possibly enlarged base space. 

In Section~\ref{subsec:Bayes} we relate our general existence result from 
Section~\ref{sec:existence} to the construction of an important class of prior distribution
families in nonparametric Bayesian inference.

For Sections~\ref{subsec:FI} and~\ref{subsec:MSE} the starting point is the observation
that a representation such as~\eqref{eq:mixdef} relates
the parametric families $\{W_\theta:\, \theta\in \Theta\}$ of probability distributions on $(E,\cE)$,
where $W_\theta$ is given by $W_\theta(\{i\})=w_\theta(i)$ for all $i\in E$, and the 
family $\{P_\theta:\, \theta\in\Theta\}$. On general
grounds the passage from the first to the second family entails an information loss. 
This can be formalized by the respective Fisher information. We obtain an integral expression 
for the classical case with Poisson distributions and noncentral chi-squared distributions with
one degree of freedom, see Proposition~\ref{prop:suffchi}. 
It turns out that at least half of the information is lost, and that this bound is asymptotically
tight as $\theta\to\infty$. Moreover, it follows that the method of moments estimator for the 
noncentrality 
parameter in the chi-squared case (which, together with some of its variants, has been 
considered in several papers) is not efficient. Further, in the general setup, 
the existence of a sufficient statistic $T$, on a possibly enlarged base
space, also leads to a comparison of experiments by conditioning on~$T$. This is carried out 
in Proposition~\ref{prop:RBunif} for the representation obtained in Section~\ref{sec:existence}
for a family of uniform distributions and the method of moments estimator. 

In Section~\ref{subsec:aseff} we note that a discrete mixture representation leads to
an embedding of the parametric family into a nonparametric one, meaning that the original
parameter $\theta$ is replaced by the probability simplex on $(E,\cE)$.
In our final result, Theorem~\ref{thm:aseff}, we show, again in the classical situation,
that the method of moments estimator for the mean functional is then asymptotically efficient 
at a large class of distributions, including the noncentral chi-squared with one degree of freedom. 
General aspects of
and comments on nonparametric maximum likelihood estimation for such classes, including the use of the EM algorithm, are collected in Section~\ref{subsec:NPMLE}.  

Section~\ref{sec:summ} concludes our work and also mentions some directions for future
research. 
Finally, an appendix contains the proofs of our results.

\section{Existence of discrete mixture representations}\label{sec:existence}
We construct a discrete mixture representation for a specific distribution family
and then use this result to obtain such representations for other families. 
We begin the first step by outlining a general approach; 
see also Section~\ref{subsec:Bayes} for a similar treatment of a problem in 
nonparametric Bayesian inference. 

Suppose that the basic space $(S,\cS)$ is the increasing limit of a 
sequence of finite spaces in the sense that $\cS$ is generated by the union 
of $\cF_n$, $n\in\bN$, where $(\cF_n)_{n\in\bN}$ is a
filtration consisting of finite $\sigma$-fields. Then each $\cF_n$ is generated by
a finite measurable partition $F_{n,1},\ldots,F_{n,k_n}$ of $S$, and these are nested.
For a prospective dominating probability measure $\nu$ we then put
\begin{equation*}
	E(\nu)\, := \, \{(n,k):\, n\in\bN, \, k=1,\ldots,k_n, \, \nu(F_{n,k})>0\},
\end{equation*}
and for $i=(n,k)\in E(\nu)$ we let $Q_i$ be the distribution with $\nu$-density
$g_{n,k} \, := \, \nu(F_{n,k})^{-1} \, 1_{F_{n,k}}$. (Here and below $1_A$ 
denotes the indicator function
of the set $A$.) For an arbitrary probability measure $\mu$
on $(S,\cS)$ with $\nu$-density $f$ we then obtain an increasing sequence of
subprobabilities via their densities $f_n=\sum_{k=1}^{k_n}a_{n,k} g_{n,k}$, with 
$a_{n,k} :=  \sup\{t\ge 0:\, \nu(tg_{n,k}\le f)=1\}$ for all $(n,k)\in E(\nu)$.
The desired discrete mixture representation then appears through the corresponding differences
if $\lim_{n\to\infty}\int(f-f_n)\, d\nu =0$. 

It should be clear that this approach can also be used for two-sided filtrations 
$(\cF_n)_{n\in\bZ}$ and countably infinite  partitions. We now carry this out for a 
specific parametric family, using the decomposition of the real line into dyadic intervals. 
Let
\begin{equation*}
\Qbin=\{k2^m:\, k\in\bZ\text{ odd}, \, m\in\bZ\}
\end{equation*} 
be the set of binary rational numbers. By convention, the notation $(a,b)$ may refer to a pair 
of real numbers or to an open interval, but the meaning should always be clear from the context.

Let $\unif(a,b)$, $-\infty<a<b<\infty$, be the uniform distribution on the interval $(a,b)$ and let
\begin{equation*}
    \Theta:=\{(a,b)\in\bR^2:\, -\infty < a < b<\infty\}.
\end{equation*}
We define a countable set $E\subset\bR^2$ by
\begin{equation}\label{eq:pairs}
    E\, :=\, \bigl\{ \bigl(k2^m,(k+1)2^m\bigr):\,k,m\in \bZ\bigr\}.
\end{equation}
Our first result now shows that the parametric family of uniform distributions
on bounded intervals of real numbers has a discrete mixture representation. 
We give a constructive proof (see the appendix), which will prove instructive later on 
when we use the representation provided by Theorem~\ref{thm:existunif}.
The proof is based on a decomposition of finite real intervals into those 
with binary rational endpoints and length an integer power of $2$. The condition in 
the theorem ensures that the decomposition is unique, but see also 
Theorem~\ref{thm:unifconv}\,(a).

\begin{theorem}\label{thm:existunif}
For $(a,b)\in \Theta$ let 
$C(a,b)$ be the set of pairs $(p,q)=(k2^m,(k+1)2^m)\in E$
with $(p,q)\subset (a,b)$ and the property that 
$(k+2)2^m>b$ if $k$ is even or $(k-1)2^m<a$ if $k$ is odd.
Then, for each $(a,b)\in\Theta$,
\begin{equation}\label{eq:reprunif0}
	\unif(a,b) \, = \, \sum_{(p,q)\in E} w_{(a,b)}(p,q) \, \unif(p,q),
\end{equation}   	
where the mixing coefficients are given by 
\begin{equation}\label{eq:mixunif} 
     w_{(a,b)}(p,q)=\begin{cases}
                       \displaystyle\frac{2^m}{b-a}, &\text{ if } (p,q)=\bigl(k2^m,(k+1)2^m\bigr)\in C(a,b),\\
                          \quad  0, &\text{ otherwise.}
                    \end{cases}
\end{equation} 
\end{theorem}

\begin{remark}\label{rem:unif}
A version of the theorem for the subfamily $\{\unif(0,\theta):\, 0<\theta < 1\}$
is of separate interest; it is also easier to state: Each $\theta$ has a unique binary 
expansion $\theta=\sum_{k=1}^{K(\theta)} 2^{-j_k(\theta)}$, 
with $K(\theta)<\infty$ if $\theta\in\Qbin$, and then,
ignoring the dependence on $\theta$ in the notation, 
\begin{equation}\label{eq:unif_onesided}
    \unif(0,\theta) \;=\;  \sum_{k=1}^{K} 2^{j_k}\, \unif\bigl(a_{k-1},a_k\bigr),
\end{equation}
with $a_0:=0$ and $a_k:=\sum_{l=1}^k 2^{-j_l}$ for $k>0$. This will be taken up in Section~\ref{subsec:MSE}.
\end{remark}

Let $(\Omega,\cA)$ and $(\Omega',\cA')$ be measurable spaces and suppose that 
$T:\Omega\to\Omega'$ is $(\cA,\cA')$-measurable.
We recall that the push-forward $P^T$ of a probability measure $P$ on 
$(\Omega,\cA)$ under $T$ is the probability measure on $(\Omega',\cA')$
given by $P^T(A)=P(T^{-1}(A))$, $A\in\cA'$. This is also known as the image of $P$
under $T$.

\begin{remark}\label{rem:structural}
The following structural property of the above discrete mixture representation 
is worth noting; it also plays a role in the last part of the proof:
A set $\{P_\theta:\, \theta\in\Theta\}$ of distributions on the Borel subsets 
of the real line is a location-scale family if the push-forward of any element under 
an affine-linear transformation $x\mapsto cx+d$, $c\not= 0$, is again an element of the family.
Clearly, this holds for the family in Theorem~\ref{thm:existunif}. The (countable) subset 
of mixing distributions that we obtain if the interval bounds are restricted
to binary rationals still enjoys this invariance property, provided that we restrict 
the shift $d$
to $\Qbin$ and $|c|$ to an integer power of $2$. 
\end{remark}

It should be clear that 
a family of distributions has a discrete mixture representation if it can written as the 
union of a finite or countably infinite family of families with this property, or as a subset. 
Below we repeatedly use two further properties: For the first of these, suppose 
that~\eqref{eq:mixdef} holds and that $\tau$ is  a probability measure on 
$(\Theta,\cB(\Theta))$. Then the $\tau$-mixture $R_\tau$ with base family 
$\{P_\theta:\, \theta\in E\}$ can be written as 
\begin{equation}\label{eq:mixmix}
     R_\tau(A):= \int P_\theta(A)\, \tau(d\theta)\; 
                   = \; \sum_{i\in E} \tilde w_\tau(i)\, Q_i(A)\quad\text{for all }A\in\cS,
\end{equation}
with $\tilde w_\tau(i):= \int w_\theta(i)\, \tau(d\theta)$ (here we implicitly assume that 
the functions $\theta\mapsto w_\theta(i)$, $i\in E$, are measurable). 
Hence a family of mixtures of the
original family again has a discrete mixture representation, even with the same family of mixing
distributions. For the second property let $(S',\cS')$ be another measurable space,
let $T:\Omega\to\Omega'$ be measurable, and assume that~$\{P_\theta:\, \theta\in\Theta\}$
has a  mixture representation $P_\theta =\int Q(y,\cdot)\, \mu_\theta(dy)$ as 
in~\eqref{eq:mixgen}. Then the following representation of the push-forwards holds,
\begin{equation}\label{eq:mixpush}
     P^T_\theta = \int Q^T(y,\cdot)\; \mu_\theta(dy)\quad\text{for all }\theta\in\Theta,
\end{equation}
where $Q^T(y,\cdot)$ denotes the push-forward of $Q(y,\cdot)$ under $T$. Clearly, the first of these
properties can be extended  to non-discrete base families, and the second can easily be specialized
to the case of discrete $E$.

As mentioned in Section~\ref{sec:intro}, noncentral distributions are of particular interest.
For these, we start with a distribution $\mu$ on the real line and write 
$P_\theta$ for the distribution of $(X+\theta)^2$ or $|X+\theta|$, where the random
variable $X$ is supposed to have distribution $\mu$. From~\eqref{eq:mixpush} it is 
clear that for the existence of a discrete mixture 
representation it is irrelevant which of these possibilities we choose.
Also, if $\mu$ is symmetric then we may assume that~$\theta\ge 0$. Below, unless 
otherwise specified, densities refer to densities with respect to the Lebesgue 
measure.

\begin{theorem}\label{thm:noncentral1}
	Suppose that $\mu$ is symmetric and has a density $f$ that is weakly decreasing 
	on $\bR_+$. Then the associated noncentral family has a discrete mixture representation.  
\end{theorem}

The proof uses a representation of $\mu$ as a mixture 
of uniform distributions on the intervals $(-y,y)$, $y>0$.
    
Both the Cauchy distribution and the double exponential distribution satisfy the assumptions 
in Theorem~\ref{thm:noncentral1} which means that, answering a question raised in Section~\ref{sec:intro},
both noncentral families have a discrete mixture representation.
In the following example we give some details for the double exponential case.

\begin{example}\label{ex:dexp}
The (standard) double exponential distribution $\Dexp(1)$ is given by its density 
$x\mapsto e^{-|x|}/2$, $x\in\bR$. It is easily checked that the corresponding
representing measure $\mu$ in the proof of Theorem~\ref{thm:noncentral1} is equal 
to the gamma distribution $\Gamma(2,1)$ with 
parameters $2$ and $1$, which has density function $x\mapsto xe^{-x}$, $x\ge 0$, so that
\begin{equation}\label{eq:dexp}
	\Dexp(1)\, = \, \int \unif(-y,y)\, \Gamma(2,1)(dy).
\end{equation}
Now let $P_\theta$ be the distribution of $|X+\theta|$, where $X\sim \Dexp(1)$ and 
$\theta\ge 0$. Then~\eqref{eq:dexp} implies 
\begin{equation}\label{eq:dexp0}
   	P_\theta \, =\, \int \unif(-y,y)^{T_\theta}\, \Gamma(2,1)(dy),
\end{equation} 
with $T_\theta(x)=|x+\theta|$. It is easy to check that the push-forwards in~\eqref{eq:dexp0}
can all be written as a uniform distribution or as a mixture of two uniform distributions. In both 
cases Theorem~\ref{thm:existunif} is applicable. 
\end{example}

In Theorem~\ref{thm:unifconv} below we consider the mixture family associated with the mixing 
distributions that appear in~\eqref{eq:reprunif0} in more detail.

\section{Representations with continuity properties}\label{sec:continuity}

We recall that noncentral chi-squared distributions 
$\chi^2_n(\theta^2)$ with $n$ degrees of freedom may be written as
\begin{equation}\label{eq:chinrepr}
     \chi_n^2(\theta^2)\, =\, \sum_{k=0}^\infty \Po(\theta^2/2)(\{k\})\, \chi^2_{2k+n},  
     \quad\theta\in\bR,
\end{equation}
where $\Po(\lambda)$ abbreviates the Poisson distribution with parameter $\lambda$ and
$\chi^2_m$ is the (central) chi-squared distribution with $m$ degrees of freedom,
i.e.\ the distribution of $X_1^2+\cdots+ X_m^2$ if~$X_1,\ldots,X_m$ are independent
standard normal random variables. This provides a discrete mixture representation for
the noncentral chi-squared distributions $\chi_n^2(\theta^2)$ with arbitrary $(n,\theta^2)\in\bN\times (0,\infty)$, and taking $n=1$ leads to such a representation for the subfamily addressed 
in the introduction. In both cases we may take $E$ to be the set of integers $i>0$ and 
$\{\chi^2_i:\, i\in E\}$ as the family of mixture distributions.

In contrast to these and the similarly explicit representations obtained in~\cite{BaGr5}
our results in Section~\ref{sec:existence} with uniforms as mixing distributions seem to
be less `usable'. In particular, they differ with respect to smoothness properties.
For example,  we may view
$\theta\mapsto (i\mapsto w_\theta(i))$ as a function on $\Theta$ with values in   
the Banach space $(\ell_1,\|\cdot\|_1)$, 
\begin{equation}\label{eq:ell}
\ell_1=\ell_1(E):=
      \Bigl\{ (a_i)_{i\in E}\in\bR^E:\, \|a\|_1:= \sum_{i\in E} |a_i|<\infty\Bigr\}.
\end{equation}
It is easy to see that the functions given by the mixing coefficients in~\eqref{eq:chinrepr} 
are  continuous. Notice though that in Theorem~\ref{thm:existunif}
and the results built on it we do not require the mixing coefficients to be continuous.  

What happens if we impose continuity assumptions on the discrete mixture representation?
We suggest to formalize these by regarding~\eqref{eq:mixdef} as taking place in a certain 
Banach space $(\bB,\|\cdot\|)$; indeed, we have already done so when passing 
from~\eqref{eq:mixdef} to~\eqref{eq:mixdef0}. In the classical case~\eqref{eq:chinrepr}
with $n>1$
we may for example use $\bB=C_b(\bR_+)$, the space of bounded continuous functions 
$f:\bR_+\to \bR$, with the supremum norm $\|f\|_\infty=\sup_{x\ge 0} |f(x)|$. The continuous
densities of $\chi^2_n$, $n>1$, are all bounded and it is easy to see that the respective 
norms tend to $0$ as $n\to\infty$, which implies boundedness of the whole family 
in~$(\bB,\|\cdot\|)$. For~the subfamily with $n=1$  we may similarly use $C_b([\epsilon,\infty])$
with an arbitrary $\epsilon>0$ (the continuous density of $\chi^2_1$ is not bounded).

The following simple result can be used to show that for a given family 
$\{P_\theta:\, \theta\in\Theta\}$
a representation of the type~\eqref{eq:mixdef} is only possible if the corresponding 
smoothness assumptions are not too strong. In it, we regard $\theta\mapsto w_\theta$ as a function
on $\Theta$ with values in the Banach space $(\ell_1,\|\cdot\|_1)$.

\begin{proposition}\label{prop:cont}
	Suppose that the representation~\eqref{eq:mixdef} holds in some Banach space $(\bB,\|\cdot\|)$, and that
	\begin{align*}
		&\theta\mapsto w_\theta \  \text{is  continuous,}\tag{C} \\
		&\{Q_i:\, i\in E\} \  \text{is bounded in  $(\bB,\|\cdot\|)$}.\tag{B}
	\end{align*}
Then, the function $\theta\mapsto P_\theta$ on $\Theta$ with values in $\bB$ is  continuous.
\end{proposition}

\begin{remark}\label{rem:cont}
(a)	Variations of this result are easily obtained. If (C) is amplified to Lipschitz continuity,
for example, then $\theta\mapsto P_\theta$ is Lipschitz continuous too. Further, if
$E=\bN$ then exponential coefficients can be introduced and (B) can be relaxed 
or amplified to the boundedness of $\rho^k Q_k$, $k\in\bN$, with some $\rho>0$
if a corresponding bound is assumed to hold for the mixing coefficients.

(b) Continuity of course also depends on the topology chosen on $\Theta$. In fact,
from a probabilistic point of view, especially in connection with the standard
model for an infinitely repeated toss of a fair coin, one might argue
for $\Theta=\{0,1\}^\infty$ instead of $\Theta=(0,1)$ in the context of 
the family $\{ \unif(0,\theta):\, 0<\theta <1\}$.
With the discrete topology on $\{0,1\}$ and 
the product topology on the new $\Theta$ the function  
$\theta\mapsto (i\mapsto w_\theta(i))$ then turns out to be  continuous; see also
Proposition~\ref{prop:RBunif}\,(b). 
\end{remark}

We revisit two noncentral families under this continuity perspective. 
Below we write $\Leb$ for the Lebesgue measure on the Borel subsets $\cB$ of $\bR$
and $\Leb_+$ for its restriction to the Borel subsets $\cB_+$ of $\bR_+$ or $\bR_+\setminus\{0\}$.
Continuity of a mixture representation refers to a notion of convergence for probability
measures.  In~\eqref{eq:mixdef} convergence of the series is the convergence of 
real numbers. As the mixing coefficients are nonnegative and summable, this  automatically 
amplifies to convergence in total variation norm of the partial sums if we regard these as 
measures. For a dominated sequence this in turn is equivalent to  $L_1$-convergence of 
the respective densities. In the special case $(S,\cS)=(\bR,\cB)$ and with dominating
measure $\Leb$ the distance of probability measures then refers to the distance of
densities, that is, of functions $f:\bR\to \bR_+$. For these, other notions, stronger
than $L^1$-convergence, can used, for example the distance based on the essential supremum,
or distances that use smoothness properties of the functions. 

\begin{example}\label{ex:irreg1}
We consider the uniform distributions $\unif(\theta,\theta+1)$, $\theta\in\bR$.
Clearly, all $P_\theta$ can be interpreted, via their densities, as elements of
\begin{equation*}
\bB \,:=\, L_\infty(\bR,\cB, \Leb)
\, = \, \bigl\{f:\bR\to \bR:\,   \|f\|_{\text{ess sup}} := \inf\{a:\, \Leb(\abs{f}>a)=0\} <\infty\bigr\}.
\end{equation*}
Further,  $\theta \mapsto P_\theta$ is $\|\cdot\|_{\text{ess sup}}$-bounded.  
For a bounded set of mixing
distributions and with continuity of the mixing coefficients, a representation of the form~\eqref{eq:mixdef} would imply that 
$\theta \mapsto P_\theta$ is continuous, which is obviously not the case: For example,
$\| P_{\theta+1/n}-P_{\theta}\|_{\text{ess sup}} \ge 1$ for all $\theta\in\bR$,  $n\in\bN$.	
\end{example}

\begin{example}\label{ex:irreg2}
Let $\mu=\Dexp(1)$ as in Example~\ref{ex:dexp}. We consider the corresponding noncentral 
distributions $P_\theta$, $\theta\ge 0$, where $P_\theta$ is the distribution of $|X+\theta|$ 
and $X$ has distribution $\mu$. Then a continuous density of $P_\theta$ is given by
\begin{align*}
	f_\theta(y)\ &= \  \frac{1}{2}\bigl(e^{-|y-\theta|}+ e^{-y-\theta}\bigr) \\
	                     &= \ \begin{cases}
	                                        e^{-\theta}(e^y+e^{-y})/2, &0<y\le\theta,\\
	                                        e^{-y}(e^\theta+e^{-\theta})/2, &\theta<y<\infty.          
	                     \end{cases}
\end{align*}
This implies that $f_\theta(y)=f_\theta(0)+\int_0^y g_\theta(z)\, dz$ for all $y\ge 0$, 
with
\begin{equation*}
g_\theta(z)\  = \ \begin{cases}    e^{-\theta}(e^z-e^{-z})/2, &0<z \le\theta,\\
                                 -e^{-z}(e^\theta+e^{-\theta})/2, &\theta<z<\infty.          
                                 \end{cases}
\end{equation*}
In particular, $f_\theta$ is differentiable in $y\not=\theta$, and the derivative has a jump of 
size 1 at $y=\theta$. We can interpret all $P_\theta,$ via the $g_\theta,$  
as elements of the Banach space $\bB \,=\, L_\infty(\bR_+,\cB_+, \Leb_+)$.
Repeating the step from~\eqref{eq:mixdef} to~\eqref{eq:mixdef0} we may
regard a representation as taking place in $\bB$. In this space, $\theta\to P_\theta$
is not continuous, so a discrete mixture representation for the noncentral family
associated with $\Dexp(1)$ satisfying the assumptions~(C) and~(B) does not exist. 
\end{example}

\section{Geometric aspects}\label{sec:geo}

We briefly sketch Dynkin's approach~\cite{Dynkin78} to convex measurable structures 
in the context of the present
situation. Let $\bM_1=\bM_1(S,\cS)$ be the set of all probability measures on $(S,\cS)$ and let
$\cB(\bM_1)$ be the $\sigma$-field on $\bM_1$ generated by the projections $P\mapsto P(A)$, $A\in\cS$. 
Then each probability measure $\Xi$ on $(\bM_1,\cB(\bM_1))$ defines a probability measure 
$P=\Psi(\Xi)$ on $(S,\cS)$, the \emph{barycenter} of $\Xi$, via
\begin{equation}\label{eq:bary}
   P(A) = \int R(A) \; \Xi(dR) \quad\text{for all } A\in\cS.
\end{equation}
Note that we integrate with respect to a measure $\Xi$ on a set of probability measures,
which means that the integration variable $R$ is a probability measure; also, 
$R(A)$ is short-hand for the function $R\mapsto R(A)$ with $A$ fixed.
We say that a subset $M$ of $\bM_1$ is \emph{barycenter convex} if 
$\Psi(\Xi)\in M$ for all probability measures $\Xi$ on $(M,\cB(M))$, where $\cB(M)$ 
denotes the trace of $\cB(\bM_1)$ on $M$. The classical notion appears if we
restrict this to measures $\Xi$ that are concentrated on a finite number of points in $M$. 
The set of all probability measures on $(\bN,\cP(\bN))$ with finite support provides an 
example of a family that is classically convex but not barycenter convex; here and 
below $\cP(A)$ denotes the power set associated with a set $A$, i.e.\ the set of all 
subsets of $A$.

Now let $(E,\cE)$ be another measurable space and let $Q$ be a transition probability 
from $(E,\cE)$ to $(S,\cS)$.  We write $\Mix\{Q(y,\cdot):\, y\in E\}$ for the set of
probability measures on $(S,\cS)$ that arise as $Q$-mixtures, see~\eqref{eq:mixgen}.
It is easy to see that $\Phi:E\to \bM_1$, $y\mapsto Q(y,\cdot)$, is 
$(\cE,\cB(\bM_1))$-measurable, and that the mixed distribution $P_\tau$ is the 
barycenter of the push-forward $\Xi=\tau^\Phi$ of $\tau$ under $\Phi$. Further, 
from the behavior of mixtures under push-forwards, see~\eqref{eq:mixpush},
it follows that such mixture families are barycenter convex; indeed,
$\Mix\{Q(y,\cdot):\, y\in E\}$ may be seen as the \emph{barycentric convex hull} of the
family $\{Q(y,\cdot):\, y\in E\}$.

The basic equation~\eqref{eq:mixdef} then says that 
$\{P_\theta:\, \theta\in\Theta\}$ is a  subset of $\Mix\{Q_i:\, i\in E\}$,
where we have written $Q_i$ instead of $Q(i,\cdot)$.
By the mixture-of-mixtures formula~\eqref{eq:mixmix} this implies
\begin{equation}\label{eq:mixsubsets}
\Mix\{P_\theta:\, \theta\in\Theta\} \subset \Mix\{Q_i:\, i\in E\},
\end{equation}
and it seems natural to call a discrete mixture representation~\eqref{eq:mixdef}  
\emph{minimal}, if
\begin{equation*}
\Mix\{P_\theta:\, \theta\in\Theta\} = \Mix\{Q_i:\, i\in E\}.
\end{equation*}
In this context, a description of the respective extreme points is of interest. 
Suppose that $\mu=\sum_{i\in E} p_iQ_i\in\Mix\{Q_i:\, i\in E\}$ is a `true' mixture
in the sense that $0<p_j<1$ for some $j\in E$.  Then $\mu$ can written as
\begin{equation*}
\mu = (1-p_j)\mu_1 + p_j\mu_2 \ \text{with }
\mu_1:= \sum_{i\not = j} \frac{p_i}{1-p_j} Q_i,\ \mu_2:= Q_j,
\end{equation*}
hence the set of extreme elements of $\Mix\{Q_i:\, i\in E\}$ is a subset of 
$\{Q_i:\, i\in E\}$. Notice that this argument uses barycenter convexity; indeed,
for sets that are convex in this stronger sense 
the classical notion of extreme points (not a non-trivial finite 
affine combination of other points of the set) and the barycentric version (not 
representable in the sense of~\eqref{eq:bary} by some $\Xi$ that is not concentrated at one 
point)  are the same.

For noncentral chi-squared distributions we have the following result.

\begin{theorem}\label{thm:simplexchi}
\emph{(a)} The representation~\eqref{eq:chinrepr} for the subfamily 
$\{\chi_1^2(\theta):\,  \theta\ge 0\}$ is not minimal, i.e.
\begin{equation*}
\Mix\{\chi_1^2(\theta):\, \theta\ge 0\} \subsetneqq \Mix\{\chi_{1+2k}^2:\, k\in\bN_0\}.
\end{equation*}
Further,  each $\chi_1^2(\eta)$, $\eta\ge 0$,  is extremal in 
$\Mix\{\chi_1^2(\theta):\, \theta\ge 0\}$,
and a minimal discrete mixture representation of $\{\chi_1^2(\theta):\, \theta\ge 0\}$ does not
exist.

\emph{(b)} The representation~\eqref{eq:chinrepr} of the family 
$\{\chi_k^2(\theta): \, \theta\ge 0, \, k\in\bN\}$ is minimal, i.e.
\begin{equation*}
\Mix\{\chi_k^2(\theta):\, \theta\ge 0,\, k\in\bN\} = \Mix\{\chi_{k}^2:\, k\in\bN\}.
\end{equation*}

\emph{(c)} If $k,l\in\bN$ have different parities, i.e.\ if $|k-l|$ is odd, then 
\begin{equation*}
\Mix\{\chi_k^2(\theta):\, \theta\ge 0\}  \cap \Mix\{\chi_l^2(\theta):\, \theta\ge 0\} =\emptyset.
\end{equation*}
\end{theorem}

The proof of the first part implies that the family	
$\Mix\{\chi_1^2(\theta):\, \theta\ge 0\}$
is a simplex, where again the notion refers to general barycenters.

We next consider the situation in Theorem~\ref{thm:existunif}, with
$E = \{(k2^m,(k+1)2^m):\, k,m\in\bZ\}$
and
\begin{equation*}
Q_{(p,q)} =\unif(p,q)\ \text{ for all } (p,q)=(k2^m,(k+1)2^m)\in E.
\end{equation*}
Again, absolutely continuity refers to the Lebesgue measure $\Leb$.

\begin{theorem}\label{thm:unifconv} 
\emph{(a)} The set $\Mix\{Q_{(p,q)}:\, (p,q)\in E\}$ has no extreme elements.

\vspace{.5mm}

\emph{(b)} Let $\mu$ be a probability measure on $(\bR,\cB)$. 
If $\mu$ is absolutely continuous with a density that is Riemann integrable
on all compact intervals, then $\mu$ is an element of $\Mix\{Q_{(p,q)}:\, (p,q)\in E\}$.

\vspace{.5mm}

\emph{(c)} There  exists a probability measure on $(\bR,\cB)$ that is absolutely 
continuous and that is not an element of 
$\Mix\{Q_{(p,q)}:\, (p,q)\in E\}$. 
\end{theorem}

The above approach can be used to obtain similar results for other families of
distributions. For example, in $\Mix\{\unif(0,\theta):\, \theta > 0\}$ the mixing distributions
are extreme, and the family contains the family $\cF$ of all distributions on $(0,\infty)$ 
with a weakly decreasing density. Obviously, $\unif(0,\theta)\in\cF$. Taken together
this shows that $\cF$, which has a discrete mixture representation by
Theorem~\ref{thm:noncentral1}, does not have a minimal discrete mixture representation.

With respect to the general approach in this section we point out that,
in contrast to the classical functional-analytic results on convexity and associated 
representations, see e.g.~\cite{Phelps}, we have not used any topological concepts
(other than those for the real line that are inherent in the Lebesgue integral). 
Nevertheless, it is interesting to compare mixture families to the 
closed convex hull of the mixing
distributions. This set of course depends on the topology in use. For example, with 
$\Mix\{\unif(0,\theta):\, \theta > 0, \, \theta\in\bQ\}$ and total variation norm (or,
equivalently, the $L_1$-norm for the respective densities), the closed convex hull
is strictly larger than the mixture family itself. Part (a) of Theorem~\ref{thm:unifconv}
may also be of interest in this connection; see also the corresponding remarks in 
the introduction.

Finally, we briefly consider what happens if we drop the assumption that 
the mixing coefficients are nonnegative.
For example, what is the closure of the linear span of the shifted Cauchy distributions 
in the space $L_1$?
A famous result from functional analysis, see e.g.~Theorem~9.5 in~\cite{RudinFA}, 
states that this set is the whole of $L_1$ if the characteristic function (Fourier transform 
of the density) has no zeroes. In the Cauchy case, this is the function $x\mapsto e^{-|x|}$, 
so that the condition is satisfied. Clearly, the corresponding mixture family is much smaller.

\section{Statistical aspects}\label{sec:stat}
The structural decomposition given by a mixture representation is closely related to the 
statistical concept of sufficiency. As in~\eqref{eq:mixgen} let $(E,\cE)$ and $(S,\cS)$ be
measurable spaces and let $Q$ be a Markov kernel from $(E,\cE)$ to $(S,\cS)$. For any probability 
measure $\tau$ on $(E,\cE)$ we may define the probability 
measure $\tau\otimes Q$ on the product space $(E \times S, \cE \otimes \cS)$ by
\begin{equation}\label{eq:product}
       \tau\otimes Q(A \times B) = \int_A Q(y,B)\,\tau(dy)
                         \quad\text{ for all }A\in \cE, \,B \in \cS.
\end{equation}
Let $T$ and $X$ be the projections $(y,x)\mapsto y$ and $(y,x)\mapsto x$ and suppose that
$D$ is a non-empty subset of the set $\bM_1(E,\cE)$ of probability measures on $(E,\cE)$.
By construction,  $T$ is then 
sufficient for the set $\{\tau\otimes Q:\, \tau\in D\}$, and
$\{P_\tau:\, \tau\in D\}$ with $P_\tau=(\tau\otimes Q)^X$ is the corresponding 
set of distributions for the second component $X$, 
with $P_\tau$ as in~\eqref{eq:mixgen}. This shows that, given a mixture representation, 
it is possible  to enlarge
the sample space such that a sufficient statistic appears. Of course, this is also quite
evident from the two-stage interpretation of mixture experiments. 
Also, we noted in the remarks following~\eqref{eq:mixgen} that a mixture representation 
always exists if we take $Q(y,\cdot)=\delta_y$. Here this corresponds to the 
statement that $X$ itself is a sufficient statistic.

Conversely, suppose that we have a set of probability measures
$\{P_\tau:\tau\in D\}$ on $(S,\cS)$, indexed by $\tau\in D\subset \bM_1(E,\cE)\}$,
with the property that on an enlarged space
$(E\times S,\cE\otimes\cS)$ there is a family $\{R_\tau:\, \tau\in D\}$ of probability measures 
such that for each $\tau\in D$ the push-forward $R_\tau^X$ of the projection $X$ on the second
coordinate is $P_\tau$, the push-forward $R_\tau^T$ of the projection $T$ on the first
coordinate is $\tau$, and $T$ is sufficient for $\{R_\tau:\tau\in D\}$.
Then if $(S,\cS)$ is a Borel space, there exists a Markov
kernel $Q$ from $(E,\cE)$ to $(S,\cS)$ such that 
\eqref{eq:mixgen} holds for all $\tau\in D;$ for more details see \cite{LeCam96}.  

In applications, $D$ is often a suitably parametrized set of distributions, say
$D=\{\tau_\theta:\theta\in\Theta\}$. Then putting $P_\theta=P_{\tau_\theta}$
the mixture representation of the parametric family
$\{P_\theta:\theta\in\Theta\}$ is obtained in the usual parametric form
\begin{equation*}
  P_\theta(\cdot)=\int Q(y,\cdot)\,\tau_\theta(dy),\quad\theta\in\Theta.
\end{equation*}

From this point of view we may rephrase our quest for a discrete mixture 
representation as the search for a discrete sufficient statistic on a possibly enlarged 
base space. Such an enlargement may indeed be necessary: If $T:S\to E$ is sufficient 
for a family $\{P_\theta:\, \theta\in\Theta\}$ of distributions on $(S,\cS)$ dominated by
some $\sigma$-finite measure~$\nu$, and if $E$ is countable, then the Neyman criterion implies
that for $\theta_1,\theta_2\in \Theta$ 
with $\theta_1\not=\theta_2$  the ratio $x \mapsto p_{\theta_1}(x)/p_{\theta_2}(x)$ 
of the associated $\nu$-densities
has only countably many values. In the classical case, with $P_\theta=\chi_1^2(2\theta)$
for $\theta > 0$, this ratio is a continuous function, and we would obtain a contradiction
to the intermediate value theorem.

A discrete mixture representation such as~\eqref{eq:mixdef} connects the statistical experiments
$\fE=(E,\cE,\{W_\theta:\, \theta \in \Theta\})$, with $W_\theta$ again given by $W_\theta(\{i\})=w_\theta(i)$ for all $i\in E$,
and $\fF=(S,\cS,\{P_\theta:\, \theta \in \Theta\})$
and also leads to an embedding of the parametric experiment $\fF$ into the nonparametric
experiment
\begin{equation*}   
   \fM=\Bigl(S,\cS,\Bigl\{\sum_{i\in E} p_i Q_i:\, p_i\ge 0\text{ for all }i\in E, \, \sum_{i\in E} p_i=1\Bigr\}\Bigr).
\end{equation*} 
Further, with the notation introduced above, we obtain the experiment
\begin{equation*}   
   \fR=\bigl(E\times S,\cE\otimes \cS,\{W_\theta\otimes Q:\, \theta\in \Theta\}\bigr).
\end{equation*} 
After a short subsection on the connection to nonparametric Bayesian inference we
investigate various statistical aspects that relate the experiments $\fE$, $\fF$, $\fM$ 
and $\fR$. In the context of comparing experiments the  connection between geometry 
and statistics can be
seen as an underlying thread, a connection that currently receives much attention
under the heading `information geometry'.  The classical case is the
use of the Hilbert space Pythagorean theorem in the proof of the 
Cram\'{e}r-Rao lower bound, which plays a role in 
Subsection~\ref{subsec:MSE}. Also, concepts from differential geometry are increasingly used
in statistics, with~\cite{Amari} 
and~\cite{Efron} being two early influential contributions. 
A specific aspect of this influence is the interpretation of Fisher information as
curvature. In the situation considered here this leads to an interpretation of the results in
Subsection~\ref{subsec:FI}
as a reduction of curvature (`flattening') in the transition from $\fE$ to~$\fF$.
Finally, convexity considerations in the context of sufficiency are basic themes 
in~\cite{Lauritzen} and~\cite{Dynkin78}.

\subsection{Priors on sets of probability measures}\label{subsec:Bayes}
We refer the reader to the  recent textbook~\cite{Ghosal-vdVaart} for a 
general introduction to nonparametric Bayesian inference. 

As in Section~\ref{sec:geo} let  $M$ be a subset of the set
$\bM_1(S,\cS)$  of all probability measures on 
$(S,\cS)$ and let $\cM$ be the $\sigma$-field on $M$ generated by the maps
$P\mapsto P(A)$, $A\in\cS$. The data $x\in S$ are regarded as a realization of
a random variable $X$,  and $M$ is the set of potential distributions for $X$.
Our prior knowledge is formalized by a probability $\Xi$  on $(M,\cM)$, 
the prior distribution. An important problem of nonparametric Bayesian inference 
is the construction of a set of such measures $\Xi$ that is flexible in the sense
that the transition to posterior distributions is manageable and does not leave the set, 
and that the Bernstein-von Mises theorem applies. 
In his seminal paper, Ferguson~\cite{Ferguson} puts special 
emphasis on Dirichlet processes for the case that $M$ is the set of all probability 
measures on $(S,\cS)$. These priors satisfy the above requirements, but
they are concentrated on the set of discrete distributions. 

Hoff~\cite{Hoff} points out the relevance of mixture representations 
for the construction of prior probabilities on convex sets of probability measures.
There, a necessary first step is the construction of a suitable mixture representation;
see~\cite[Section 4]{Hoff} for an interesting variety of worked examples.
Here we start with a mixture representation. For example, in
the statistical experiment $\fM$ we have
\begin{equation*}   
M\; = \; \Bigl\{\sum_{i\in E} p_i Q_i:\, p_i\ge 0\text{ for all }i\in E, \, \sum_{i\in E} p_i=1\Bigr\},
\end{equation*} 
and the construction of probability measures on $M$ is essentially equivalent to the
construction of probability measures on $\bM_1(E,\cP(E))$. 
In the classical situation and in many other cases of interest, $E=\bN_0$, and the latter
problem can be approached via a stick-breaking procedure.
The detour via a discrete mixture representation means that the posterior
distribution obtained from the prior and the data  would still be absolutely 
continuous with respect to any measure dominating $M$ (in fact, even 
smoothness properties of the densities would be retained).

An important alternative to Dirichlet processes are tree-based priors. 
To see the connection
to discrete mixture representations we assume, as at the beginning of
Section~\ref{sec:existence}, that $\cS$ is generated by the union 
of $\sigma$-fields  $\cF_n$, $n\in\bN$, where $(\cF_n)_{n\in\bN}$ is a
filtration and each $\cF_n$, $n\in\bN$,
is generated by a finite partition $\{F_{n,1},\ldots,F_{n,k_n}\}$. This structure leads
to a tree with node set $(n,k)$, $n\in\bN$, $k=1,\ldots,k_n$, and edges between
$(n,k)$ and $(n+1,j)$ whenever $F_{n,k}\supset F_{n+1,j}$. A suitable assignment
of probabilities $p_{(n,k),(n+1,j)}$ to these edges then defines a probability $P$
on $(S,\cS)$, with $ P\bigl(E_{n+1,j}\big| E_{n,k}\bigr) \, = \, p_{(n,k),(n+1,j)}$ for all edges $\{(n,k),(n+1,j)\}$. Choosing random assignments leads to tree-based priors $\Xi$,
such as tail-free processes and P\'{o}lya trees. 

Conversely, choosing the weights via the conditional probabilities provides a 
tree-based representation for a given distribution $P$. The approximation $P_n$
that results if the conditional probabilities  of the nodes are multiplied up to depth 
$n$ is simply  the restriction of $P$ to  the $\sigma$-field $\cF_n$. 
For discrete mixture representation, however, we need a 
sequence of approximations that is increasing. Such a sequence of subprobability 
measures can be obtained via a dominating measure, as explained at the beginning of
Section~\ref{sec:existence}.

\subsection{Fisher information}\label{subsec:FI}
We first consider the experiments $\fE$ and $\fF$ defined above, where the parameter set
$\Theta$ is assumed to be an open subset of $\bR^d.$ In what follows we confine ourselves
to giving the definition of the Fisher information matrix for the experiment $\fF.$ Let $f(\cdot,\theta)$
be the density of the distribution $P_\theta$ with respect to some $\sigma$-finite measure
on $(S,\cS).$ We denote by $\theta_1,\ldots,\theta_d$ the components of the column vector $\theta\in\Theta$.
Under suitable standard regularity conditions, see, e.g. \cite{Lehmann},
the integrals
\begin{equation}\label{eq:FIclassical}
i_{\fF;jk}(\theta)=\int \frac{\partial}{\partial \theta_j}\log f(x,\theta)
                          \ \frac{\partial}{\partial \theta_k}\log f(x,\theta)
                          \,P_\theta(dx),\quad 1\le j,k \le d,
\end{equation}                        
exist and it is $i_\fF(\theta)=\left (i_{\fF;jk}(\theta)\right )_{1\le j,k\le d}$
a symmetric positive definite $d\times d$ matrix which is called the
Fisher information matrix associated with the experiment $\fF.$ 
Assuming that the corresponding standard conditions hold for the experiment $\fE$,  
there is a Fisher information matrix $i_\fE(\theta)$ associated with $\fE$ as well.
The experiments $\fE$ and $\fF$ are related in the sense that
there is Markov kernel $Q$ from $(E,\cE)$ to $(S,\cS)$
such that~\eqref{eq:mixdef} holds. Following Le~Cam~\cite{LeCam96}, the experiment $\fF$ 
is then said to be reproducible from the experiment $\fE$ (or, equivalently, $\fE$ is said to be better than $\fF$).
Notice  that the standard terminology that $\fE$ is sufficient for $\fF$ is not adopted by Le Cam.
Under certain regularity conditions it then follows that
\begin{equation}\label{Fisher:matrix}
  i_\fE(\theta) - i_\fF(\theta)~~\text{is positive semidefinite for each}~\theta\in\Theta; 
\end{equation}
see \cite{GoelDeGroot79} or \cite{Stone61}. The latter author deals with the
one-dimensional case $d=1.$ Then \eqref{Fisher:matrix} reduces to
\begin{equation*}
  i_\fE(\theta) \ge  i_\fF(\theta)\quad \text{for each}~\theta\in\Theta. 
\end{equation*}
It is of interest to quantify $i_\fE(\theta) - i_\fF(\theta),$ which can be viewed as the
loss of information when switching 
from $\fE$ to $\fF$. We tackle this problem in the special case arising with
the mixture representation
\begin{equation}\label{eq:classical}
   \chi^2_1(2\theta) \, =\, \sum_{k=0}^\infty \Po(\theta)(\{k\})\, \chi^2_{2k+1},\quad\theta\in (0,\infty).
\end{equation}
The above conditions for~\eqref{eq:FIclassical} are then satisfied for the 
experiments $\fE$ and $\fF$ that are specified in the following result. 

\begin{proposition}\label{prop:suffchi}
In the special case of the experiments
\begin{equation*}
 \fE=\bigl(\bN_0,\cP(\bN_0),\{\Po(\theta):\theta \in (0,\infty)\}\bigr)
\end{equation*}
and
\begin{equation*}
   \fF= \bigl(\bR_{>0},\cB_{>0},\{\chi_1^2(2\theta):\theta \in (0,\infty)\}\bigr)
\end{equation*}
it holds that
\begin{equation}\label{eq:FIintegral}
    i_\fF(\theta) = \frac{1}{2\theta}\int_{0}^{\infty}x\,\tanh^2\bigl((2\theta x)^{1/2}\bigr)\,
             g_{1,\theta}(x)\,dx\; -\;  1.
\end{equation}
In particular, with $r(\theta):=\theta i_\fF(\theta)$ for all $\theta>0$,
\begin{equation}\label{eq:prop_r}
  \lim_{\theta\to 0} r(\theta)=0,\quad \lim_{\theta\to\infty}r(\theta)=\frac{1}{2},
\end{equation}
and $r$ is strictly increasing. Finally,
\begin{equation}\label{eq:FIinequ}
          i_\fF(\theta)\;\le\; \frac{1}{2}i_\fE(\theta)\ \text{ for all }\theta >0.
\end{equation}
\end{proposition}

It follows from~\eqref{eq:FIinequ} that the transition from  
$\fE$ to $\fF$ 
causes a loss of information of more than fifty percent, and~\eqref{eq:prop_r} implies 
that  this bound is asymptotically tight as $\theta\to\infty$.

For a general discussion of estimation procedures for estimating the unknown parameter $\theta$
based on observations of $X$ we refer to \cite{Kubokawa17}.

\subsection{MSE reduction}\label{subsec:MSE}
\begin{figure}
	\includegraphics[width=0.7\linewidth]{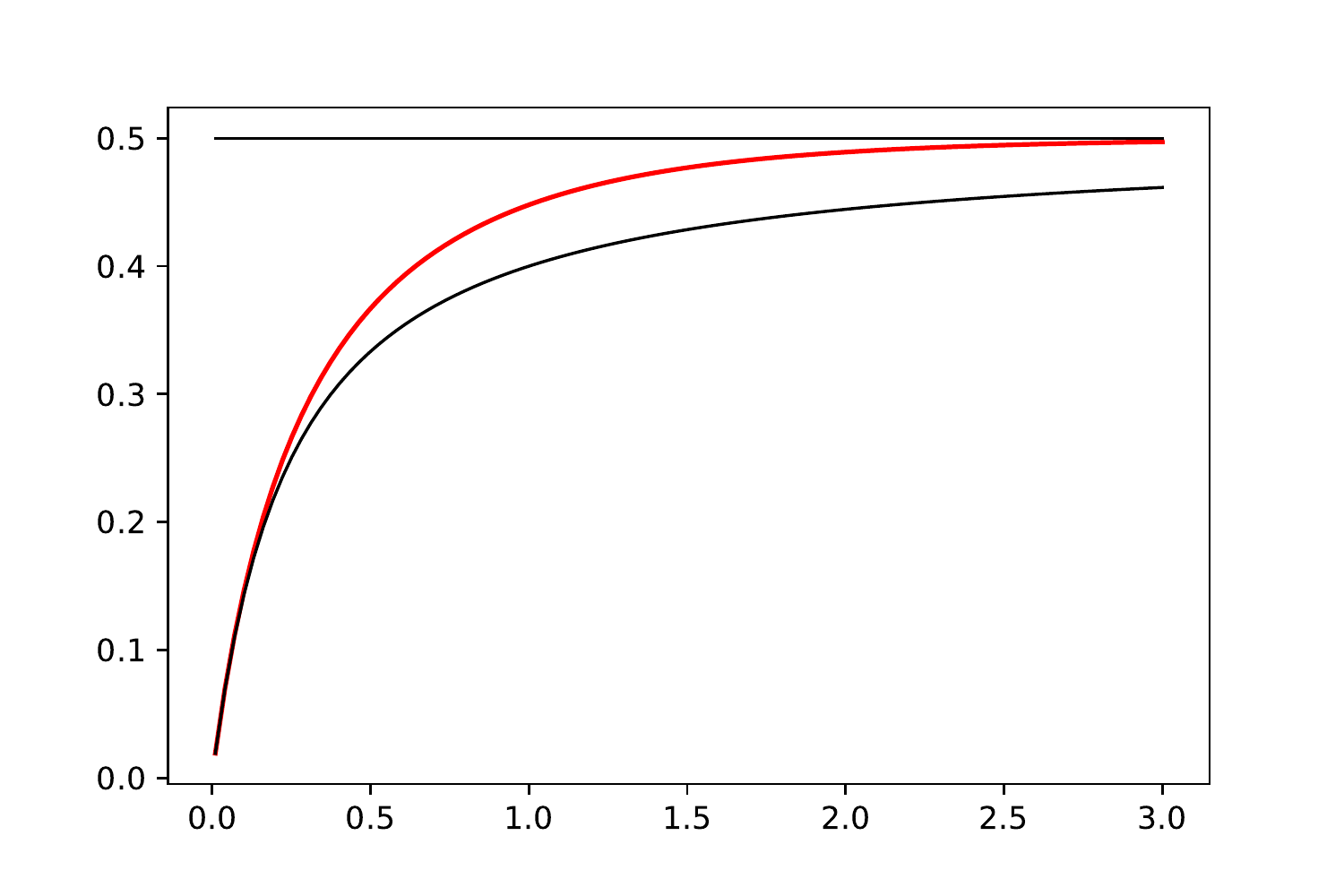}
	\caption{The function $\theta\mapsto\theta i_{\fF}(\theta)$ (red) and its bounds, 
		see~\eqref{eq:FIbounds}.}
	\label{fig:suffchi}
\end{figure}

A quantitative comparison of $\fF$ and $\fR$, the original and 
the enlarged experiment, can be obtained through the reduction of the mean squared 
error (MSE) of estimators by conditioning on the sufficient 
statistic, a procedure sometimes called `Rao-Blackwellization'.
In the classical case~\eqref{eq:classical}, with $\Theta=(0,\infty)$ and $P_\theta = \chi_1^2(2\theta)$, 
the enlargement leads to the product space $(\bN_0 \times \bR_{>0}, \cP(\bN_0) \otimes \cB_{>0})$,
with the distributions $R_\theta$ specified by
$R_\theta(\{k\} \times A) = \Po(\theta)(\{k\})\cdot \chi_{2k+1}^2(A)$ for all $k\in \bN_0$, $A \in \cB_{>0}$.
Again, $T$ and $X$ are the canonical projections on $\bN_0$ and $\bR_{>0}$ respectively,
and  $T$ is sufficient for $\{R_\theta :\, \theta\in\Theta \}$. In this situation, $\hat\theta:=(X-1)/2$
is an unbiased estimator for the unknown parameter $\theta$, and conditioning on the sufficient statistic 
leads to the estimator $\tilde\theta = E_\theta[\hat\theta|T]=T$. The MSE reduction may
be given explicitly in this situation:
\begin{equation*}
\MSE_\theta(\tilde\theta)\,=\, \theta\ <\ \frac{1}{2}+ 2\theta \,=\, \MSE_\theta(\hat\theta).
\end{equation*}
We note that, in view of the Cram\'{e}r-Rao lower bound for the variance of an unbiased
estimator, the formula for $\MSE_\theta(\hat\theta)$ can also be used to obtain a lower 
bound for the Fisher information $i_\fF$ considered in the previous subsection. In fact,
inspecting the geometric argument in the proof of the bound it is not difficult to show that
the inequality is strict, so that the result of the previous subsection may be augmented as
follows,
\begin{equation}\label{eq:FIbounds}
	\frac{1}{2+1/(2\theta)}\, < \, \theta i_{\fF}(\theta)\, < \, \frac{1}{2}\quad\text{for all }\theta>0,
\end{equation} 
see also Figure~\ref{fig:suffchi} where the integral in~\eqref{eq:FIintegral} 
has been computed numerically.
Note that \eqref{eq:FIbounds} leads to an alternative proof of the second limit relation stated
in~\eqref{eq:prop_r}. Further,
standard arguments show  that $\var(\hat\theta)>i_{\fF}(\theta)^{-1}$ implies that 
the estimator $\hat\theta$ is not asymptotically efficient 
and thus inferior to the maximum likelihood estimator; 
see also Subsection~\ref{subsec:aseff}.

As another, somewhat less straightforward example we consider the family of uniform distributions
$P_\theta=\unif(0,\theta)$, $\theta\in\Theta:=(0,1)$, on $(S,\cS)=((0,1),\cB_{(0,1)})$ 
together with the mixture representation in Remark~\ref{rem:unif}.
In order to carry out Rao-Blackwellization for the moment estimator $\hat\theta=2X$, 
$X\sim P_\theta$, we use the binary representation: For a given $\theta\in\Theta$, 
let $K(\theta)$, $j_k(\theta)$ and $a_k(\theta)$ be as in~\eqref{eq:unif_onesided}. 
Let $E:=\Theta\cap\Qbin$, and for $t\in E$
let $d(t):=2^{-K(t)}$ where $K(t)$ is obtained from the binary representation of $t$. 
Given $\theta\in \Theta$
we assign the weight $w_\theta(a_m(\theta))=d(a_m(\theta))/\theta$ to each 
$a_m(\theta)$, $m=1,\ldots,K(\theta)$, 
and put $w_\theta(t)=0$ otherwise. Further, with
each $t\in E$ we associate the mixing distribution $Q(t,\cdot)=\unif(t-d(t),t)$.
Then the discrete mixture representation \eqref{eq:unif_onesided} may be written as 
\begin{equation*}
  P_\theta=\unif(0,\theta)=\sum_{t\in E} w_\theta(t)\, \unif\left (t-d(t),t\right ).
\end{equation*}
For example, if $\theta=3/4$ then  $K(\theta)=2$, $a_1(\theta)=1/2$, $a_2(\theta)=3/4$,
$d(a_1(\theta))=1/2$, $d(a_2(\theta))=1/4$, $w_\theta(1/2)=2/3$, $w_\theta(3/4)=1/3$, 
and we arrive at 
\begin{equation*}
    \unif\Bigl(0,\frac{3}{4}\Bigr)\; 
             =\; \frac{2}{3}\,\unif\Bigl(0,\frac{1}{2}\Bigr)\ +\  \frac{1}{3}\,\unif\Bigl(\frac{1}{2},\frac{3}{4}\Bigr).
\end{equation*}
On the product space $(E \times (0,1), \cP(E) \otimes \cB_{(0,1)})$
we obtain the probability measures $R_\theta$ specified by
\begin{equation*}
    R_\theta(\{t\} \times A) = 
         w_\theta(t)\,\unif\left (t-d(t),t\right )(A),\ t\in E,\,A \in \cB_{(0,1)}.
\end{equation*}
With $T$ and $X$ the canonical projections on $E$ and $(0,1)$ respectively, 
$T$ is sufficient for $\{R_\theta : \theta\in\Theta \}$. For the moment estimator 
$\hat\theta=2X$ conditioning on the sufficient statistic now leads to the estimator 
$\tilde \theta:= E_\theta[\hat\theta|T]=2T-d(T)$. As in Proposition~\ref{prop:suffchi}
we first give a general formula and then derive some properties of the function of interest.

\begin{proposition}\label{prop:RBunif}
\emph{(a)} With the notation as introduced above, the  mean squared
error of the conditioned moment estimator $\tilde\theta=2T-d(T)$ is given by 
\begin{equation}\label{eq:MSEtilde}
	\MSE_\theta(\tilde\theta) = \var_\theta(\tilde{\theta})
	\, =\, \frac{1}{\theta}  
             \sum_{m=1}^{K(\theta)} a_m(\theta) a_{m-1}(\theta)\bigl(a_m(\theta)-a_{m-1}(\theta)\bigr).
\end{equation}
\emph{(b)} The function $\theta\mapsto \phi(\theta):=\var_\theta(\tilde{\theta})$ is continuous on $\Theta\setminus E$, and on $E$ it is right continuous and has left hand limits.
Moreover, 
\begin{equation}\label{eq:sprung}
     \phi(\theta)-\phi(\theta-)\; =\;  -\frac{1}{7 q}\, 2^{-2L+1}
\end{equation} 
if $\theta$ is a binary rational of the form $\theta=q2^{-L}$, with $L,q\in\bN$ and $q$ odd.
\end{proposition}

The property `right continuous, with left hand limits' is often abbreviated to `c\`adl\`ag'.

\begin{figure}
	\includegraphics[width=0.7\linewidth]{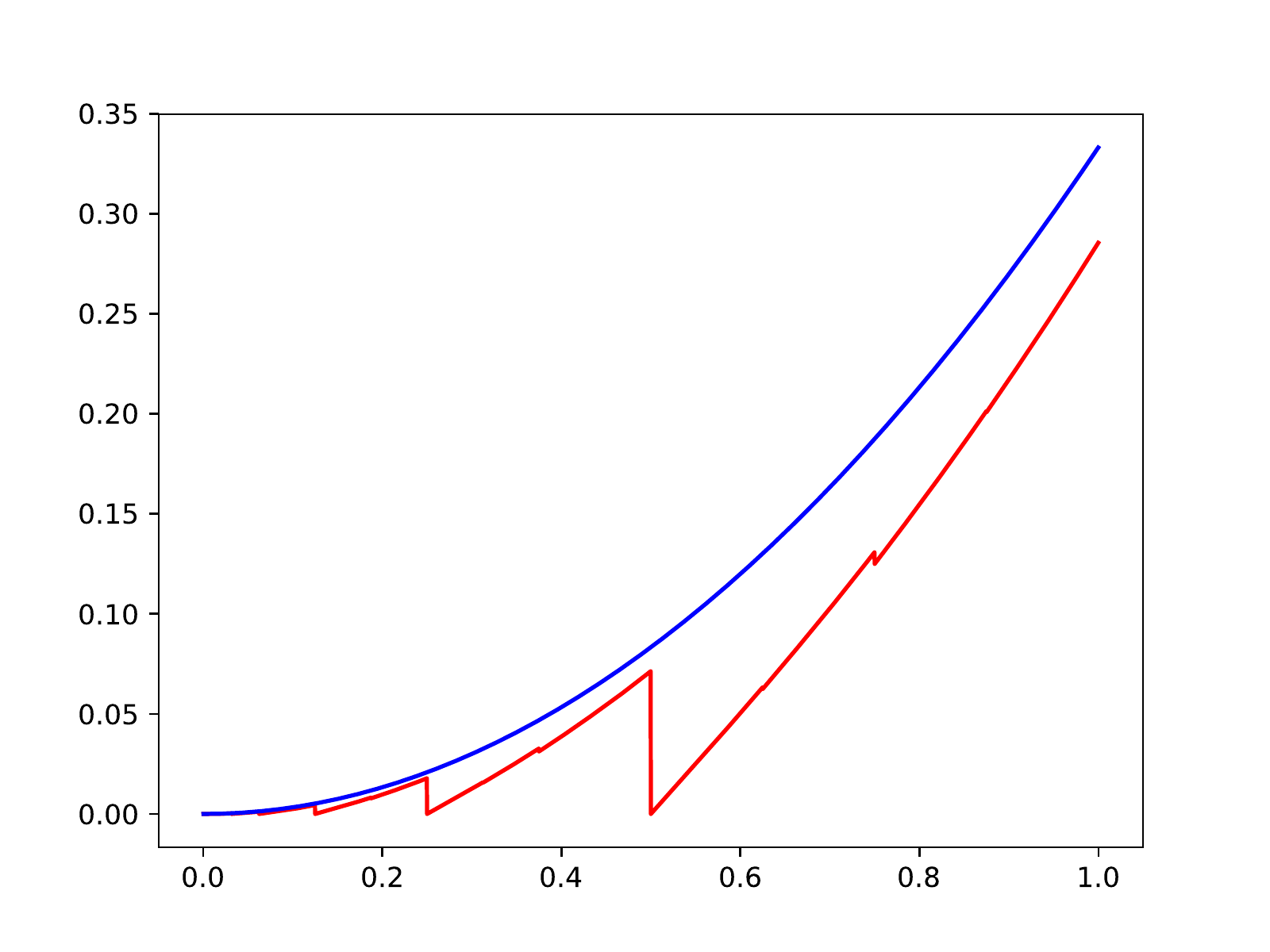}
	\caption{Mean squared errors associated with $\hat\theta$ (blue) and $\tilde\theta$ (red)}
	\label{fig:suffunif}
\end{figure}

For the unconditioned estimator $\hat\theta$ we have $\var_\theta(\hat\theta)=\theta^2/3$.
Both estimators are unbiased. Figure~\ref{fig:suffunif} shows the variances (mean squared errors) of  
$\hat\theta$ and $\tilde\theta$ for $\theta$ an integer multiple of $2^{-12}$, with the interpolation
justified by Proposition~\ref{prop:RBunif}\,(b). The figure 
suggests a self-similarity property: Indeed, as $X\sim P_{\theta/2}$ and $2X\sim P_\theta$	
are equivalent, it follows that both variance functions have the scaling property 
$\phi(\theta/2)= \phi(\theta)/4$. Further, if $\theta=2^{-m}$ for some $m\in\bN$ then the
variance is zero, as the distribution of $T$ is degenerate for such values of the parameter.

\subsection{An asymptotically efficient estimator for the mean functional}
\label{subsec:aseff}

We consider the classical case~\eqref{eq:classical}. In the corresponding 
mixture model we then have
\begin{equation*}
	\Theta =\Bigl \{p=(p_k)_{k\in\bN_0}:\,  p_k\ge 0\text{ for all }k\in\bN_0,\,
	\sum_{k=0}^\infty p_k=1\Bigr\},
\end{equation*} 
and the statistical experiment $\fM$ is based on the family
\begin{equation*}
    \cP\, =\, \Bigl\{P_p:\, P_p=\sum_{k=0}^\infty p_k \chi_{2k+1}^2,\, p\in\Theta\Bigr\}
\end{equation*}    
of mixture distributions 
parametrized by $p\in\Theta$. The family  $\left \{\chi_1^2(2\theta):\,\theta>0\right \}$
underlying $\fF$ is a subfamily of $\cP$. 

Our aim is to show that the moment estimator for the mean functional, which we know
from the previous subsections to not be efficient in the context of $\fF$, is asymptotically 
efficient at each $P_p\in\cP$ that has non-vanishing mixing coefficients with finite variance. 
We refer the reader to~\cite[Sections 1.2, 2.1 and 3.1]{Vaart} for an exposition of the 
general theory needed here. 
In comparison with the situation in Proposition~\ref{prop:suffchi} there are three main
differences: First, instead of a parametric family with a parameter $\theta$
specifying the distribution, we now have a one-dimensional parameter function $\kappa:\cP\to \bR$,
where $\kappa(P)$ does not fully characterize $P$. Secondly, instead of the variance
of an estimator $\hat\theta$ for $\theta$ we now consider the variance of the limiting 
normal distribution in an associated central limit theorem for a sequence of estimators
for $\kappa(P)$.
The third point requires some machinery. The basic idea is to consider dominated 
one-dimensional submodels $\{P_t:\, 0\le t < \epsilon\} \subset \cP$ with $P_0=P$ and
then use the derivative at $t=0$, in analogy to~\eqref{eq:FIclassical}. This leads to a tangent
set, where the geometry is that of a Hilbert space of square integrable functions.  
Again speaking somewhat informally, a lower bound 
can then be obtained from the supremum of the associated second moments. Finally, an estimator 
sequence is said to be asymptotically efficient at $P$ for the parameter function $\kappa$ if asymptotic variance and lower bound are the same for this $P$.

Let $\cP_0$ be the subset of distributions $P_p\in \cP$ with  $p=(p_k)_{k\in\bN_0}$ 
satisfying the moment condition $\sum_{k=0}^\infty k^2p_k<\infty$ and also the condition
that $p_k>0$ for each $k\in\bN_0.$
Consider the mean functional $\kappa:\cP_0\rightarrow \bR$ defined by
\begin{equation*}
   \kappa(P_p)=\sum_{k=0}^\infty k p_k,\,P_p\in \cP_0.
\end{equation*}
Let $X_1,\dots,X_n,\dots$ be a sequence of independent and identically
distributed random variables with distribution $P_p\in \cP_0$ which is assumed to be unknown. 
Then, generalizing the moment estimator that already appeared in Subsection~\ref{subsec:MSE}, 
\begin{equation*}
	T_n:=\frac{1}{2}\bigl(\overline{X}_n-1\bigr)\, =\,\frac{1}{2}\Bigl(\frac{1}{n}\sum_{i=1}^n X_i-1\Bigr)
\end{equation*}
leads to a sequence $(T_n)_{n\in\bN}$ of unbiased estimators for $\kappa(P_p)$, 
and the central limit theorem shows that, as $n\to\infty$,
\begin{equation}\label{eq:unbiasedreg}
  \sqrt{n}\left (T_n-\kappa(P_p)\right )
                      \xrightarrow{\cD} N\Bigl(0,v(P_p)+\kappa(P_p)+\frac{1}{2}\Bigr),
\end{equation}
where $v(P_p)=\sum_{k=0}^\infty \left (k-\kappa(P_p)\right )^2p_k$. 

\begin{theorem}\label{thm:aseff}
  For estimating $\kappa(P_p),$ the estimator sequence $(T_n)_{n=1}^\infty$ is asymptotically
  efficient at each $P_p\in\cP_0$.
\end{theorem}

If $P_p=\chi_1^2(2\theta)$  then $v(P_p)=\kappa(P_p)=\theta$ so that 
the variance of the limit distribution of $\sqrt{n}\left (T_n-\kappa(P_p)\right )$
is equal to $2\theta+\frac{1}{2}$, in accordance with the value found in 
Subsection~\ref{subsec:MSE}.

\subsection{NPMLE and EM}\label{subsec:NPMLE}

As we pointed out in the introduction a distribution $P_\theta$ with the mixture 
representation $P_\theta=\sum_{i\in E} w_\theta(i)\, Q_i$
may be seen as the distribution of the outcome $X$ in a two-stage experiment: First, choose
$Y$ in $E$ with distribution (mass function) $w_\theta$, then, if $Y=i$, choose $X$ according 
to $Q_i$. In a sample $x_1,\ldots,x_n$ of size $n$ from $P_\theta$ we may then think of the
corresponding $y_1,\ldots,y_n$ as hidden (or latent) variables, or as missing covariates, 
see~\cite[Sect.~1.3.5 and 1.3.6]{Lindsay}. A non-parametric generalization leads to the 
problem of estimating the mixing probabilities $p_i$, $i\in E$, in a representation
$P=\sum_{i\in E} p_i\, Q_i$ of the unknown distribution $P$, with the $Q_i$ known, 
from the data $x_1,\ldots,x_n$. We may of course regard the sequence $p=(p_i)_{i\in E}$ as
a parameter, where the parameter space is now the probability simplex on $E$, as in the
previous subsection.

The EM algorithm, see~\cite{DLR} and~\cite[Sect. 3.4]{Lindsay} can be used to 
obtain a sequence $p(l)$, $l=1,2,\ldots\,$ of approximations to the corresponding 
non-parametric maximum likelihood estimator (NPMLE). We give the details for the 
classical case, where $E=\bN_0$ and $Q_k$ is the central chi-squared distribution with 
$2k+1$ degrees of freedom, with continuous density $g_k$.
The log-likelihood function $L_X$ is then given by
\begin{equation}\label{eq:llX}
    L_X(p)\, =\, \sum_{i=1}^{n} \log\Bigl(\sum_{k=0}^{\infty} p_k\, g_k(x_i)  \Bigr).
\end{equation} 
With the corresponding $y$-values known, the log-likelihood function would be
\begin{equation}\label{eq:llXY}
    L_{X,Y}(p)\, =\, \sum_{i=1}^{n} \log\bigl(p_{y_i}\, g_{y_i}(x_i)\bigr)
              \, =\, \sum_{i=1}^{n} \log p_{y_i}\, +\, r(x_1,\ldots,x_n;y_1,\ldots,y_n),
\end{equation} 
where the function $r$ does not depend on $p$. In the EM algorithm, we obtain the next 
approximation $p(l+1)$ for the NPMLE  from the current value $p(l)$ as the argmax (M-step) 
of the conditional
expectation (E-step) of $L_{X,Y}(p)$ given $x_1,\ldots,x_n$. 

The E-step boils down to the calculation of $E_{p(l)} \bigl[\log p_{Y_i}|X_i=x_i\bigr]$.
Under $P_{p(l)}$   the joint distribution of $X_i$ and $Y_i$ has density 
$(x,k)\mapsto p_k(l)g_k(x)$ with respect 
to the product of Lebesgue measure on $\bR_+$ and counting measure on $\bN_0$, so that 
the conditional probability of $Y_i=k$ given $X_i=x$ becomes
\begin{equation*}
h_x(k) \, =\, \frac{p_k(l)\, g_k(x)}{\sum_{j=0}^\infty p_j(l) \, g_j(x)},\quad x>0,\, k\in\bN_0.
\end{equation*}
The M-step then requires maximization of the function 
\begin{equation*}
p \, \mapsto\, \sum_{k=0}^\infty q_n(k)\, \log p_k,\ \text{ with } 
q_n(k) := \frac{1}{n}\sum_{i=1}^n h_{x_i}(k).
\end{equation*}
As $q_n$ is a probability mass function, Gibbs' inequality can be used to show that the desired argmax
is given by $p(l+1)=q_n$. Using this in~\eqref{eq:llXY} we arrive at
\begin{equation}\label{eq:EMstep}
p_k(l+1) \, =\, \frac{1}{n}\sum_{i=1}^n  \frac{p_k(l)\, g_k(x_i)}{\sum_{j=1}^\infty p_j(l) \, g_j(x_i)}
\quad \text{for all } k\in\bN_0.
\end{equation}
For chi-squared mixing distributions the factors $e^{-x/2}$ cancel in~\eqref{eq:EMstep} and, 
more importantly, the infinite sums can be truncated to the range from 0 to 
\begin{equation*}
K := K(x_1,\ldots,x_n) \, :=\, \bigl\lceil (1+ \max\{x_1,\ldots,x_n\}) / 2 \bigr\rceil.
\end{equation*}
To see this we note that
\begin{equation*}
\frac{g_{k-1}(x)}{g_{k}(x)} \; =\;   \frac{2\;\Gamma(k + 1/2)}{x\;\Gamma(k-1/2)}
\; =\, \frac{2k -1}{x} 
\quad\text{for all } x>0, \, k\in\bN,
\end{equation*}
Hence, if some $k> K$ appears in~\eqref{eq:llX}, then the corresponding mass $p_k$ can 
be shifted to the left without decreasing the value of $L_X$ . 

\begin{figure}
	\includegraphics[width=0.48\linewidth]{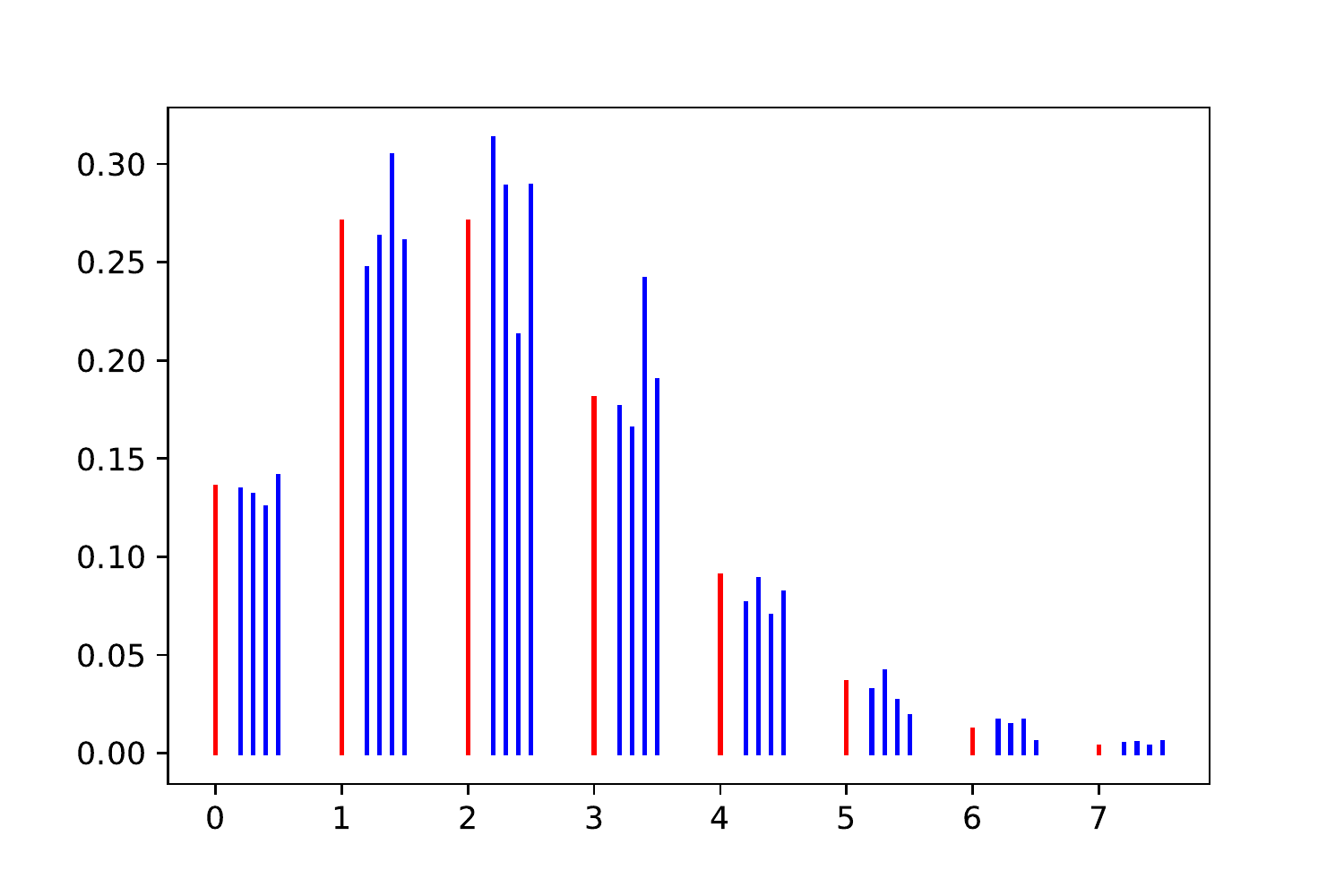}	\includegraphics[width=0.48\linewidth]{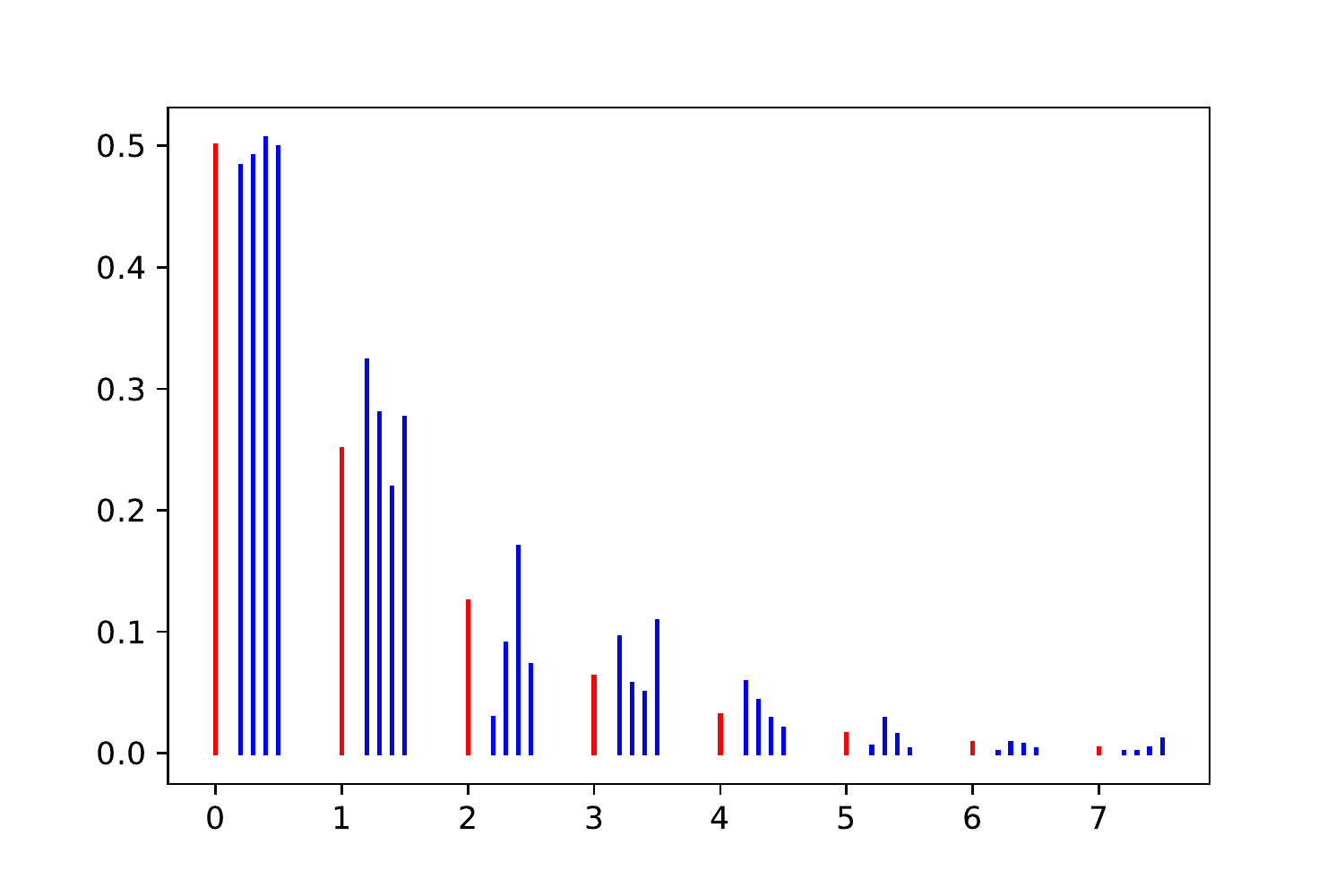}
	\caption{Simulation results (see text)}
	\label{fig:sim}
\end{figure}

The left part of Figure~\ref{fig:sim} shows the results of four 
simulation experiments, each with $n=10000$ observations, where the base distribution 
is a noncentral chi-squared with one degree of freedom and noncentrality parameter $\theta=2$.
The vertical red lines represent the probability mass $p_k$ of the true mixing distribution, 
which is Poisson with parameter~2, at the positions $k=0,1,\ldots,7$. 
The four estimates of these masses are computed numerically by using the EM algorithm 
and are shown in blue slightly to the right. 

For noncentral distribution families the nonparametric version incorporates a variant of the 
deconvolution problem~\cite[Sect.~1.3.19]{Lindsay}: In the chi-squared case, we regard the 
data as realizations of random variables $X=(Y+Z)^2$ with $Y,Z$ independent, $Y$ standard normal, and an unknown distribution $\mu$ for $Z$. The procedure outlined above can be 
applied and leads  to an estimator for the mixed Poisson distribution associated with~$\mu$. 
This in turn can then be used to obtain an estimator for $\mu$, by combining NPMLE and EM
again, for example. 

The right hand part of Figure~\ref{fig:sim} shows the results of four 
simulation experiments where $Z=2\sqrt{W}$ with $W$ exponentially distributed 
with parameter~1. The corresponding mixed Poisson distribution is the geometric
distribution on $\bN_0$ with parameter $1/2$. Again, the vertical red lines
represents the true mixing distribution, the four estimates are slightly to the right 
and in blue, and each simulation is based on $n=10000$ observations.
We should point out that a rather large number $n$ of observations is required in order to
obtain estimators with small mean squared error. 

Given a sequence $X_1,\dots,X_n,\dots$ of independent and identically
distributed random variables with distribution $P_p=\sum_{k=0}^\infty p_k\chi_{1+2k}^2,$
where the parameter $p=(p_k)_{k\in\bN_0}$ in the probability simplex on $\bN_0$ is unknown, 
the non-parametric maximum likelihood estimator $p_{\rm NPMLE,n}$ of $p$ based on
$X_1,\dots,X_n$ exists. Using the general consistency statement given in 
\cite[Theorem 5.3]{pfanzagl}
we deduce that the sequence $(p_{\rm NPMLE,n})_{n\in\bN}$ is strongly consistent for $p$.

\section{Summary and outlook}\label{sec:summ}
We have investigated geometrical and statistical aspects, and their interaction,
in the context of discrete mixture representations. It turned out that there are
connections to a variety of theoretical and applied topics, including
\begin{itemize}
\item[--] tree-based constructions of probability measures, with connections to
      nonparamnetric Bayesian inference (Theorem~\ref{thm:existunif}, 
      Section~\ref{subsec:Bayes}),
\item[--]  the view of mixture representations as curves in an infinite-dimensional
      space (Section~\ref{sec:continuity}),
\item[--]  Dynkin's barycentric approach to convexity (Section~\ref{sec:geo}),
\item[--]  the effect of mixing on Fisher information (Section~\ref{subsec:FI}) and 
      mean squared error (Section~\ref{subsec:MSE}), and
\item[--]  algorithmic aspects, notably the use of the EM algorithm (Section~\ref{subsec:NPMLE}).
\end{itemize}

We have almost exclusively considered the one-dimensional case, i.e.~families of probability
distributions on the real line, leaving multivariate situations to a future paper. Another
aspect that seems to be worth a separate investigation is the structure of the family
$\{Q_i:\, i\in E\}$ of mixing distributions in a representation such as~\eqref{eq:mixdef}.
In the classical case~\eqref{eq:chinrepr}, with noncentral chi-squared distributions,
these are the convolution powers of one specific probability measure, and we noted
in Remark~\ref{rem:structural} that the mixing distributions 
in Theorem~\ref{thm:existunif} constitute 
a location-scale family of a specific type. Taken together, the multivariate case and
the structure of the mixing family are also of interest for applications in the general
area of stochastic processes, with relations to the Dynkin isomorphism 
(see~\cite[Section 4.2]{BaGr5}). In an applied context it would certainly be interesting 
to remove the assumption in Section~\ref{subsec:NPMLE} that the mixing family is known.
This would lead to a statistical inverse problem, and structural assumptions could be used, 
if applicable, to reduce its complexity.

\vspace{.5cm}\textbf{Acknowledgment.} We thank the associate editor and the two referees
for their detailed and constructive comments, which have led to a considerable improvement
of the paper.

\section*{Appendix: Proofs}\label{sec:appendix}

\textbf{Proof of Theorem~\ref{thm:existunif}.\ } Consider the infinite rooted binary tree
where each node has one left and one right descendant. Any $x\in (0,1)$ that is not a binary
rational defines a unique infinite path through this tree, starting at the root node,
and moving to the left or right descendant if the next digit in its binary expansion is
0 or 1 respectively. With $x=1/3$, for example, we would move to the left and to the right
in odd and even steps alternatingly.
Our first aim is to obtain a decomposition of $(a,b)\subset (0,1)$ into intervals with binary 
rational endpoints. For this we label the nodes of tree by appropriately chosen
intervals, using~\eqref{eq:pairs}, and then collect these along the paths given by
$a$ and $b$.

In order to carry this out formally, let $a,b\in\bR\setminus \Qbin$, $0<a<b<1$, be given, 
with 
binary expansions $a=\sum_{k=1}^\infty a_k2^{-k}$, $b=\sum_{k=1}^\infty b_k2^{-k}$ where
$a_k,b_k\in\{0,1\}$ for all $k\in\bN$. Let $j_1<j_2<\cdots$ be the positions of the digit~1
in the expansion of $b$. The corresponding partial 
sums $b(l) := \sum_{k=1}^{j_l} b_k 2^{-k}$, $l\in\bN$, then
provide a strictly increasing sequence of binary rational numbers with limit $b$. 
By construction, $2^{j_{l+1}}\bigl(b(l+1)-b(l)\bigr)=1$,
$2^{j_{l+1}}b(l+1)-1$ is even, and $\left (2^{j_{l+1}}b(l+1)-1+2\right )2^{-j_{l+1}} = b(l+1)+2^{-j_{l+1}}>b$.

Similarly, and now writing $j_1<j_2<\cdots$ for the positions of the digit 0 in the expansion
of $a$, and with $a(l) := 2^{-j_l} + \sum_{k=1}^{j_l} a_k 2^{-k}$, $l\in\bN$,
we obtain a strictly decreasing sequence $(a(l))_{l\in\bN}\subset \Qbin$ with limit $a$ and 
the property that $2^{j_{l+1}}a(l+1)$ is odd, and $\left (2^{j_{l+1}}a(l+1)-1 \right )2^{-j_{l+1}} = a(l+1)-2^{-j_{l+1}}<a$
for all $l\in\bN$.

Let $K:=\min\{k\in\bN:\, a_k=0,\, b_k=1\}$. Clearly, $K<\infty$ as $a<b$. With 
$a(K):=\sum_{k=1}^K a_k2^{-k}$, $b(K):=\sum_{k=1}^K b_k2^{-k}$, we have $(a(K),b(K))\in C(a,b)$,
$\bigl(a(l+1),a(l)\bigr)\in C(a,b)$ for all $l$ with $j_l(a)\ge K$,
$\bigl(b(l),b(l+1)\bigr)\in C(a,b)$ for all $l$ with $j_l(b)\ge K$,
and we finally arrive at the interval decomposition
\begin{equation*}
 (a,b) = \bigl(a(K),b(K)\bigr)\;\sqcup\!\! \bigsqcup_{\{l:\, j_l(b)\ge K\}}\!\!  \bigl(b(l),b(l+1)\bigr)\;\sqcup\!\! \bigsqcup_{\{l:\, j_l(a)\ge K\}}\!\!  \bigl(a(l+1),a(l)\bigr)
                    \;\sqcup\, N,
\end{equation*}
with $N$ countable. This implies
\begin{align*}
(b-a)\,\unif(a,b)\  &=\ (b(K)-a(K))\,\unif(a(K),b(K)) \\
                  &\qquad\qquad +\ \sum_{\{l:\, j_l(b)\ge K\}}(b(l+1)-b(l))\,\unif(b(l),b(l+1)) \\
                  &\qquad\qquad +\ \sum_{\{l:\, j_l(a)\ge K\}}(a(l)-a(l+1))\,\unif(a(l+1),a(l)). 
\end{align*}

It remains to remove the restrictions on $a$ and $b$. For elements of $\Qbin$ the binary 
representation is not unique, but if we agree to choose the representation with a finite
number of digits $1$ in the case of $b$ and a finite number of digits $0$ in the case of $a$
then the only change in the above argument is that the respective sequence of $j$-values 
would be finite. 

In order to remove the assumption $0<a,b<1$ we first describe the way in 
which the representation interacts with affine transformations  $T:\bR\to\bR$ 
of the form $T(x)=2^k x+r$, with $k\in\bZ$ and $r\in\Qbin$.
Such functions are bijections, and their inverse is of the same type. 
Further, they leave $E$ invariant,
if applied to the components of the pair, and it is easy to check that
\begin{equation*}
\bigl(T(p),T(q)\bigr)\in C\bigl(T(a),T(b)\bigr)
    \quad \Longleftrightarrow \quad (p,q)\in C(a,b)
\end{equation*}
for all $(a,b)\in\Theta$, $(p,q)\in E$.
Hence~\eqref{eq:reprunif0} would imply that
\begin{equation}
    \unif\bigl(T(a),T(b)\bigr) \, = \, \sum_{(p,q)\in E} w_{(T(a),T(b))}(p,q) \, \unif(p,q),
\end{equation}
which means that the representation then also holds for $\bigl(T(a),T(b)\bigr)$.
Any pair $(a,b)\in\Theta$ can be transformed by such $T$'s into the unit interval, 
so the desired representation holds for all $(a,b)\in\Theta$.

\smallbreak
\textbf{Proof of Theorem~\ref{thm:noncentral1}.\ }
It is well-known and follows easily from Shepp's theorem~\cite[p.158]{FellerII} 
that a probability measure $P$ on the positive half-line
with weakly decreasing density $f$ (which we may take to be left continuous) may be written
as a mixture of uniform distributions on the intervals $[0,y]$, $y > 0$. For completeness we 
give a simple and direct argument: First, a $\sigma$-finite measure $\rho$ on $(\bR_+,\cB_+)$ 
can be defined via $\rho([y,\infty))=f(y)$ for all $y > 0$. Then let $\nu$ be the measure
with density $y\mapsto y$ with respect to $\rho$. Is it easy to check that $\nu$ has total mass~1.
Also, for all $x > 0$,
\begin{equation*}
   f(x) =   \int_{[x,\infty)}\frac{1}{y} \,\nu(dy)
                                   =  \int_{(0,\infty)}\frac{1}{y} 1_{[x,\infty)}(y) \,\nu(dy)
                                   = \int_{(0,\infty)}\frac{1}{y} 1_{[0,y]}(x) \,\nu(dy),
\end{equation*}
and the integrand in the last expression is a density for $\unif(0,y)$.

With $\mu$ and $f$ as in the statement of the theorem we may write 
$\mu=\frac12 (\mu_0+\mu_0^T)$ with $T(x)=-x$, and with
$\mu_0$ a distribution that is concentrated on the positive
half-line and has a weakly decreasing density there. Using the first
part, we get
\begin{equation}\label{eq:reprsymm}
\mu = \int_{(0,\infty)} \unif(-y,y)\, \nu(dy),
\end{equation}
where $\nu$ is the mixing measure for $\mu_0$. 

We now apply the general principles represented by~\eqref{eq:mixmix} 
(repeated mixtures) and~\eqref{eq:mixpush} (behaviour under push-forwards): 
First, if $X$ has distribution $\mu$ and $Q_\theta$ is the distribution of $X+\theta$, 
representation~\eqref{eq:reprsymm} implies that 
\begin{equation}\label{eq:reprshift}
 Q_\theta  = \int_{(0,\infty)} \unif(\theta-y,\theta+y)\, \nu(dy).
\end{equation}
Thus, the family 
$\{Q_\theta:\,\theta\in\bR\}$ has a discrete mixture representation in terms of the countable set
of uniform distributions on intervals with binary rational endpoints introduced in connection
with Theorem~\ref{thm:existunif}. 
The transfer argument for push-forwards, with $T(x)=|x|$ or $T(x)= x^2$, now completes the proof.  	

\smallbreak
\textbf{Proof of Proposition~\ref{prop:cont}.\ }
Assume without loss of generality that $E=\bN$. Let $C:=\sup_{k\in \bN}\|Q_k\|$. 	
Then, for all $\eta,\theta\in\Theta$,
\begin{align*}
 \|P_\theta-P_\eta\| \  
        &=\ \lim_{K\to\infty}\Bigl\|  \sum_{k=0}^K w_\theta(k) Q_k
                                           - \sum_{k=0}^K w_\eta(k) Q_k\Bigr\| \\  
    &\le \  C\,\lim_{K\to\infty} \sum_{k=0}^K \bigl|w_\theta(k)-w_\eta(k)\bigr| \\ 
       &\le \  C\; \|w_\theta-w_\eta\|_1.\qedhere
\end{align*}
\smallbreak
\textbf{Proof of Theorem~\ref{thm:simplexchi}.\ }
(a) Suppose that $\mu\in \Mix\{\chi_{1+2k}^2:\, k\in\bN_0\}$ has two mixture representations,
\begin{equation*}
\mu \, = \, \sum_{k=0}^{\infty} p_k \chi_{1+2k}^2 \text{ and } 
\mu\, = \, \sum_{k=0}^{\infty} q_k \chi_{1+2k}^2, 
\end{equation*}
with $p_k,q_k\ge 0$ and $\sum p_k=\sum q_k=1$. Then the respective densities must be 
equal almost everywhere, so that the continuous versions agree on $(0,\infty)$. This gives
\begin{equation*}
\sum_{k=0}^\infty \frac{p_k}{\Gamma\bigl(k+\frac{1}{2}\bigr)\, 
	2^{k+\frac{1}{2}}}\, x^{k-\frac{1}{2}} e^{-x/2}
\; =\; \sum_{k=0}^\infty \frac{q_k}{\Gamma\bigl(k+\frac{1}{2}\bigr)\, 
	2^{k+\frac{1}{2}}}\, x^{k-\frac{1}{2}} e^{-x/2}
\text{ for all } x>0.
\end{equation*} 
Multiplying both sides by $x^{\frac{1}{2}} e^{x/2}$ and using the uniqueness of
the coefficients in a power series representation we get $p_k=q_k$ for all $k\in\bN_0$.	
Hence the mixture representation of elements of  $\Mix\{\chi_{1+2k}^2:\, k\in\bN_0\}$ is unique.

Now suppose that $\mu$ is a mixture of noncentral chi-squared distributions with one 
degree of freedom and mixing measure $\nu$ on the parameter space $\Theta=[0,\infty)$. 
The classical representation~\eqref{eq:chinrepr} with $n=1$ then gives
\begin{equation*}
\mu=\sum_{k=0}^\infty p_k\chi^2_{1+2k} \ \text{ with }
p_k = \int e^{-\lambda} \frac{\lambda^k}{k!} \, \nu(d\lambda)\ \text{  for all } k\in \bN_0,
\end{equation*}
so that the mixing coefficients are the probabilities of a mixed Poisson distribution. In
order to prove the strict subset relation it is therefore enough to name a distribution on $\bN_0$
that is not a mixed Poisson distribution---such as $\delta_1$. In particular,
\begin{equation*}
\chi_3^2 \; \in\;  \Mix\{\chi_{1+2k}^2:\, k\in\bN_0\}\setminus 	\Mix\{\chi_1^2(\theta):\, \theta\ge 0\} .
\end{equation*}
	
To finish the proof of part (a) we use the uniqueness of the representation to show that each 
$\chi^2_1(\eta)$, $\eta>0$, is an extreme element of the convex set
$\Mix\{\chi_1^2(\theta):\, \theta\ge 0\}$. 
Indeed, suppose that $\chi_1^2(\eta) = \alpha \mu_1 + (1-\alpha)\mu_2$, and that $\nu_1,\nu_2$ 
are the parameter distributions for $\mu_1,\mu_2$. But then $\chi_1^2(\eta)$ would itself be
a mixture of  the distributions $\chi_1^2(\theta)$, $\theta\ge 0$, with mixing measure 
$\nu:= \alpha\nu_1 + (1-\alpha)\nu_2$. The uniqueness now implies that $\nu$ is concentrated 
at $\eta$, so that $\mu_1$ and $\mu_2$ are equal. 

This shows that  the set of extreme elements of $\Mix\{\chi_1^2(\theta):\, \theta\ge 0\}$ is 
uncountable, in contrast to the set of extreme elements of $\Mix\{\chi_{1+2k}^2:\, k\in\bN_0\}$.  

(b) We know from~\eqref{eq:chinrepr} that each $\chi_k^2(\theta)$, $\theta\ge 0$ and $k\in\bN$, 
has a representation in terms of central chi-squared distributions, so that, by the
mixture-of mixtures property~\eqref{eq:mixmix}, the left set is
contained in the set on the right side of the assertion. However, 
$\chi^2_k=\chi^2_k(0)$, which means that the mixing distributions are elements of the left set,
hence so are their mixtures. 
	
(c) Any element $\mu$ of both sets $\Mix\{\chi_{k+2n}^2:\, n\in\bN_0\}$ and 
$\Mix\{\chi_{l+2n}^2:\, n\in\bN_0\}$ must have two representations of the form
\begin{equation*}
\mu \, = \, \sum_{n=0}^\infty p_n\chi_{k+2n}^2 \,
\text{ and }  \mu\, = \, \sum_{n=0}^\infty q_n\chi_{l+2n}^2. 
\end{equation*}
As in the proof of (a) above, this would lead to 
\begin{equation*}
\sum_{n=0}^\infty \frac{p_n}{\Gamma\bigl(n+\frac{k}{2}\bigr)\, 
	2^{n+\frac{k}{2}}}\, x^{n+\frac{k}{2}-1} e^{-x/2}
\; =\; \sum_{n=0}^\infty \frac{q_n}{\Gamma\bigl(n+\frac{l}{2}\bigr)\, 
	2^{n+\frac{l}{2}}}\, x^{n+\frac{l}{2}-1} e^{-x/2}
\end{equation*}
for all $x>0$. Multiplication by $x^{\frac{1}{2}}e^{x/2}$ and substituting $x=y^2$ then gives
\begin{equation*}
\sum_{n=0}^\infty \frac{p_n}{\Gamma\bigl(n+\frac{k}{2}\bigr)\, 
	2^{n+\frac{k}{2}}}\, y^{2n+k -1} 
\; =\; \sum_{n=0}^\infty \frac{q_n}{\Gamma\bigl(n+\frac{l}{2}\bigr)\, 
	2^{n+\frac{l}{2}}}\, y^{2n+l-1}
\text{ for all } y>0.	
\end{equation*}
If $k$ and $l$ have different parities then the coefficients of $y$
vanish for all odd powers on one side of the equation, and for all even powers
on the other.

\smallbreak
\textbf{Proof of Theorem~\ref{thm:unifconv}.\ }
(a) On general grounds, the set of extreme elements of $\Mix\{Q_{(p,q)}:\, (p,q)\in E\}$ is a
subset of $\{Q_{(p,q)}:\, (p,q)\in E\}$. For an arbitrary  $(p,q)=(k2^m,(k+1)2^m)\in E$ we have
\begin{equation*}
\bigl(k2^{m},(k+1)2^{m}\bigr)=\bigsqcup_{r=1}^{\infty} 
                     \bigl((2^{r+1}k+1)2^{m-r-1},(2^{r+1}k+2)2^{m-r-1}\bigr)\; \sqcup \; N
\end{equation*}
with $N$ countable, so that $Q_{(p,q)}$ may be written as 
\begin{equation*}
  Q_{(p,q)}\; =\; \sum_{r=1}^\infty 2^{-r}
    \unif\bigl((2^{r+1}k+1)2^{m-r-1},(2^{r+1}k+2)2^{m-r-1}\bigr)\ 
       =\ \frac{1}{2}\,\mu_1+\frac{1}{2}\,\mu_2,
\end{equation*}
where $\mu_1:=\unif\bigl((4k+1)2^{m-2},(4k+2)2^{m-2}\big )$ and
\begin{equation*}
  \mu_2\,:=\,\sum_{r=2}^\infty 2^{-(r-1)}\unif\bigl((2^{r+1}k+1)2^{m-r-1},(2^{r+1}k+2)2^{m-r-1}\bigr)
\end{equation*}
are two different elements of $\Mix\{Q_{(p,q)}:\, (p,q)\in E\}$. 

(b) Let $f$ be a  density of $\mu$ that is Riemann integrable on all compact intervals. 
For a fixed $M\in\bZ$ consider the partition of the compact interval $I_M=[M,M+1]$
into intervals $I_{M,m,k}$, $k=1,\ldots,2^{m}$, of length $2^{-m}$. 
As these partitions are nested, the step functions $g_{M,m}$, $m\in \bN$, given by
\begin{equation*}
     g_{M,m} :=\sum_{k=1}^{2^{m}} a_{M,m,k} 1_{I_{M,k}} \quad \text{ with }
                  a_{M,m,k}:=\inf\bigl\{f(x):\, x\in I_{M,m,k}\bigr\} 
\end{equation*}
are increasing; let $g_{M,\infty}$ be their supremum. Then $0\le g_{M,\infty}\le f$,
and
\begin{equation*}
\int_{I_M} g_{M,\infty} \, d\Leb = \int_{M}^{M+1} f(x) \, dx = \int_{I_M} f \, d\Leb
\end{equation*} 
as $f$ is Riemann integrable on $I_M$. For each of the intervals $D=I_{M,m,k}$ the function 
$\Leb(D)^{-1}1_D$ is a density of $\unif(D)$, which is an element of 
$\Mix\{Q_{(p,q)}:\, (p,q)\in E\}$ by Theorem~\ref{thm:existunif}. Further, the differences 
$g_{M,m+1}-g_{M,m}$ are linear combinations with nonnegative coefficients of the indicators of 
$I_{M,m+1,k}$, $k=1,\ldots,2^{m+1}$. Taken together this shows that $g_{M,\infty}/\mu(I_M)$ is 
the density of an element of $\Mix\{Q_{(p,q)}:\, (p,q)\in E\}$ whenever $\mu(I_M)>0$. As 
$\mu$ itself is a mixture of these, an appeal to the repeated mixture 
property~\eqref{eq:mixmix} completes the proof of~(b).

(c) We start with a familiar construction from measure theory, 
see e.g.\  \cite[Problem 2.5 b), p.65]{Benedetto76}.
Let $q_k$, $k\in \bN$, be an enumeration of the rational numbers in the interval $(0,1)$.
Choose for each $k\in\bN$ some $n_k\in\bN,\,n_k> k+1,$ such that 
\begin{equation*}
I_k:=\left (q_k-2^{-n_k},q_k+2^{-n_k}\right )\subset (0,1).
\end{equation*}   
Then the union $A:=\bigcup_{k=1}^\infty I_k$ is a non-empty open subset of $[0,1]$, the
Lebesgue measure $\alpha:=\Leb(A)$ of which is 
\begin{equation*}
0<\alpha\le \sum_{k=1}^\infty \Leb(I_k)\le 2\sum_{k=1}^\infty 2^{-n_k}<\sum_{k=1}^\infty 2^{-k}=1.
\end{equation*}
The set $A$ is dense in $[0,1]$. Let $\mu$ be the uniform distribution on the complement 
$[0,1]\setminus A$ of $A$, with density $f:=(1-\alpha)^{-1}(1-1_A)$. The support of $\mu$ 
does not contain any interval of positive length, in contrast to the support of the elements of 
$\Mix\{Q_i:\, i\in E\}$. 

\smallbreak
\textbf{Proof of Proposition~\ref{prop:suffchi}.\ }
Denote by $g_{n,\theta}$ the continuous density of the $\chi_n^2(2\theta)$ 
distribution, i.e.,  
\begin{equation*}
   g_{n,\theta}(x)=\sum_{k=0}^\infty e^{-\theta}\frac{\theta^k}{k!}f_{2k+n}(x),~x>0.
\end{equation*}
Then
\begin{align*}
\frac{\rm d}{{\rm d}\theta}\log g_{1,\theta}(x)=\frac{g_{3,\theta}(x) - g_{1,\theta}(x)}{g_{1,\theta}(x)},
\end{align*}
and therefore
\begin{align*}
   i_\fF(\theta)\, 
           = \, \int_{0}^{\infty} \left (\frac{\rm d}{{\rm d}\theta}
                        \log g_{1,\theta}(x)\right )^2g_{1,\theta}(x)\,dx
       &\, = \, \int_{0}^{\infty}\left(\frac{g_{3,\theta}(x)}{g_{1,\theta}(x)}-1\right)^2
                                                    \,g_{1,\theta}(x)\,dx\\
       &\, = \, \int_{0}^{\infty}\left (\frac{g_{3,\theta}(x)}{g_{1,\theta}(x)}\right )^2
                                                      \,g_{1,\theta}(x)\,dx - 1.
\end{align*}
For $\nu\in (-1,\infty)$ let
\begin{equation*}
     I_\nu(x)=(x/2)^\nu\sum_{k=0}^\infty \frac{(x/2)^{2k}}{k!\Gamma(\nu+k+1)},\ x>0,
\end{equation*}
be the modified Bessel function of the first kind of order $\nu$. Using
\begin{align*}
   g_{1,\theta}(x) \; &=\;  \frac{1}{2}\exp(-\theta-x/2)(x/2\theta)^{-1/4}
                              I_{-1/2}\left ((2\theta x)^{1/2}\right ),\\
   g_{3,\theta}(x) \; &=\; \frac{1}{2}\exp(-\theta-x/2)(x/2\theta)^{1/4}
                              I_{1/2}\left ((2\theta x)^{1/2}\right ),
\end{align*}
see, e.g., \cite{SaxenaAlam82}, we get
\begin{equation*}
   i_\fF(\theta)\, = \, \int_{0}^{\infty}\frac{x}{2\theta}\left (\frac{I_{1/2}
             \left ((2\theta x)^{1/2}\right )}{I_{-1/2}\left ((2\theta x)^{1/2}\right )}\right )^2
                                             g_{1,\theta}(x)\,dx - 1.
\end{equation*}
Finally, from
$I_{-1/2}(x)=\left (2/(\pi x)\right )^{1/2}\cosh x$ and
$I_{1/2}(x)= \left (2/(\pi x)\right )^{1/2}\sinh x$, $x>0$, 
we deduce that  
\begin{align*}
i_\fF(\theta)&\ =\ \frac{1}{2\theta}\int_{0}^{\infty}x\,\tanh^2\left ((2\theta x)^{1/2}\right )\,
             g_{1,\theta}(x)\,dx - 1\\
           &\ <\  \frac{1}{2\theta}\int_{0}^{\infty}x\hspace*{0.5mm}g_{1,\theta}(x)\,dx - 1\quad
           =\ \frac{1}{2\theta}.
\end{align*}
The Fisher information in $\fE$ is known to be $i_\fE(\theta) = 1/\theta$ 
for $\theta > 0$. 

To prove the properties of the function $r$ we write  
\begin{align*}
  r(\theta)& \ = \ \frac{1}{2}\int_0^\infty x\tanh^2\bigl((2\theta x)^{1/2}\bigr)\,g_{1,\theta}(x)\,dx-\theta\\
           & \ = \ \frac{1}{2}\left (1-s(\theta)\right ),
\end{align*}
with
\[s(\theta)=\int_0^\infty x\left (1-\tanh^2\bigl((2\theta x)^{1/2}\bigr)\right )\,g_{1,\theta}(x)\,dx.\]
Noticing that
\begin{equation*}
    1-\tanh^2\bigl((2\theta x)^{1/2}\bigr ) \ 
       = \ 4\frac{\exp\left (-2(2\theta x)^{1/2}\right )}
              {\left (1+\exp\left (-2(2\theta x)^{1/2}\right )\right )^2} \ 
     \le \ 4\exp\left (-2(2\theta x)^{1/2}\right )
\end{equation*}
and
\begin{align*}
   g_{1,\theta}(x) \ &= \ (2\pi x)^{-1/2}e^{-\theta}e^{-x/2}\cosh\left ((2\theta x)^{1/2}\right )\\
         &\le \ (2\pi x)^{-1/2}e^{-\theta}e^{-x/2}\exp\left ((2\theta x)^{1/2}\right )
\end{align*} 
we deduce that
\begin{equation*}
   0 \ \le \  x\left (1-\tanh^2\bigl((2\theta x)^{1/2}\bigr)\right )\,g_{1,\theta}(x)
  \ \le  4(2\pi)^{-1/2}e^{-\theta}x^{1/2}e^{-x/2}.
\end{equation*}
Therefore,  
\begin{equation*}
   0 \ \le \ s(\theta)
       \le \ 4(2\pi)^{-1/2}e^{-\theta}\int_0^\infty x^{1/2}e^{-x/2}\,dx \ 
        = \ 4e^{-\theta},
\end{equation*}
so that $\lim_{\theta\to\infty}s(\theta)=0$. Additionally, due to
\begin{equation*}
   \lim_{\theta\to 0} x\left (1-\tanh^2\bigl((2\theta x)^{1/2}\bigr)\right )\,g_{1,\theta}(x) \ 
         = \ x f_1(x)\ \text{ for each }x>0,
\end{equation*}
by dominated convergence, we have that $\lim_{\theta\to 0} s(\theta) = \int_0^\infty x f_1(x)\,dx=1$.

Finally, let $0<\theta_1<\theta_2<\infty$. To show that $s(\theta_1)>s(\theta_2)$ we write
\begin{equation}\label{eq:monoton}
  s(\theta_1) \ = \ \int_0^\infty \frac{1-\tanh^2\bigl((2\theta_1 x)^{1/2}\bigr)}
                                       {1-\tanh^2\bigl((2\theta_2 x)^{1/2}\bigr )}\
                                  \frac{g_{1,\theta_1}(x)}{g_{1,\theta_2}(x)}\ 
                  x\left (1-\tanh^2\bigl((2\theta_2 x)^{1/2}\bigr)\right )\,g_{1,\theta_2}(x)\,dx.
\end{equation}
Due to 
\begin{equation*}
    \frac{1-\tanh^2\bigl((2\theta_1 x)^{1/2}\bigr)}{1-\tanh^2\bigl((2\theta_2 x)^{1/2}\bigr )}
     \ = \
    \frac{\exp\left (-2(2\theta_1x)^{1/2}\right )
          \left (1+\exp\left (-2(2\theta_2x)^{1/2}\right )\right )^2}
         {\left (1+\exp\left (-2(2\theta_1x)^{1/2}\right )\right )^2 
                     \exp\left (-2(2\theta_2x)^{1/2}\right )}
\end{equation*}
and
\begin{equation*}
   \frac{g_{1,\theta_1}(x)}{g_{1,\theta_2}(x)}
  \ = \
  \exp(\theta_2-\theta_1)\frac{\cosh\left ((2\theta_1x)^{1/2}\right )}
                              {\cosh\left ((2\theta_2x)^{1/2}\right )}
\end{equation*}
it holds that
\begin{equation*}
 \frac{1-\tanh^2\bigl((2\theta_1 x)^{1/2}\bigr)}{1-\tanh^2\bigl((2\theta_2 x)^{1/2}\bigr )}\
               \frac{g_{1,\theta_1}(x)}{g_{1,\theta_2}(x)}
               \ = \
               \exp(\theta_2-\theta_1)\frac{\cosh\left ((2\theta_2x)^{1/2}\right )}{\cosh\left ((2\theta_1x)^{1/2}\right )}
               \ > \ 1
\end{equation*} 
for each $x>0$. Thus, by \eqref{eq:monoton}, $s(\theta_1)>s(\theta_2)$.
\smallbreak

\smallbreak
\textbf{Proof of Proposition~\ref{prop:RBunif}.\ }
(a) Let $\theta\in\Theta$ be given. We first assume that $K(\theta)=\infty$ and 
suppress the dependence on $\theta$ where convenient. 
In terms of the sequences $(j_k)_{k\in\bN}$ and $(a_m)_{m\in\bN_0}$ we have 
\begin{equation*}
   R_\theta(T=a_m)	\, =\, \frac{a_m-a_{m-1}}{\theta}\, =\, \frac{1}{\theta 2^{j_m}}
                         \quad\text{for all } m\in\bN
\end{equation*}
and $R_\theta(T=t)=0$
for all other elements $t$ of $E$. Conditionally on $T=a_m$, $X$ is uniformly distributed on the interval 
$(a_{m-1},a_m]$. With these definitions we obtain
\begin{align*}
   E_\theta\tilde\theta\
      \ &=\ \sum_{m=1}^\infty R_\theta(T=a_m) E[2X|T=a_m]\
      =\ \frac{1}{\theta}\sum_{m=1}^\infty (a_{m}-a_{m-1})(a_{m}+a_{m-1})\\
     &=\ \frac{1}{\theta}\lim_{M\to\infty}a_M^2\ =\ \theta,
\end{align*} 
which, of course,  is also immediate from the tower property of conditional expectations.
Similarly, 
\begin{align*}
   \theta\, E_\theta\tilde\theta^2 \
       &=\ \sum_{m=1}^\infty (a_{m}-a_{m-1})(a_m + a_{m-1})^2\\
       &=\ \sum_{m=1}^\infty (a_m^3-a_{m-1}^3)\ 
                 + \ \sum_{m=1}^\infty a_m a_{m-1}(a_m-a_{m-1}).
\end{align*} 
From the above expression for $E_\theta\tilde\theta$ we get $\lim_{m\to\infty}a_m^3=\theta^{3}$,
which completes the proof of~\eqref{eq:MSEtilde} in the case that $K(\theta)=\infty$. For binary rational parameter values the same
arguments apply, only that the sums are now finite and no limits are needed.

\smallbreak
(b)	We first consider
the function $\Psi(\theta)= w_\theta$ that maps $\Theta$ to the set of probability measures 
on $\bN$ together with the total variation distance
or, equivalently via probability mass functions, to the space $(\ell_1(\bN), \|\cdot\|_1)$. 
Clearly, $\Psi=\Phi_1\circ\Phi_2$
where $\Phi_2$ maps $\theta=\sum_{k=1}^\infty b_k2^{-k}$ to the sequence 
$b(\theta)=(b_k)_{k\in\bN}\in\{0,1\}^\bN$ of digits in the binary expansion that has 
infinitely many 0's,  and $\Phi_1(b)(\{k\})= b_k2^{-k}/\theta$, $k\in\bN$. On the intermediate
space of 0-1 sequences we consider the product topology, where a sequence of sequences converges 
if the respective components converge, which here means that the components do not change 
from some sequence index onwards (which may depend on the index of the component; 
informally, all components eventually 'freeze').

The outer function $\Phi_1$  is continuous in view of Scheff\'e's lemma. For any 
$\theta\in \Theta\setminus E$  there are infinitely many 0's (and 1's) in $b(\theta)$.  
If $b_k(\theta)=0$ and $b_{k+1}(\theta)=1$ then $\theta$ is strictly 
inside a binary rational interval of length $2^{-k}$. Hence, if $\theta_n\to\theta$ then
$\theta_n$ is also contained in this interval from some $n_0$ onwards, which in turn implies that 
the first $k$ digits remain constant. This shows that the inner function $\Phi_2$ 
is continuous on $\Theta\setminus E$.
In a similar fashion we obtain that $\Phi_2$ is right continuous on $E$, and that in $E$ 
its left hand limits exist.

The conditioned estimator $\tilde \theta=2T-d(T)$ is bounded. Hence Lebesgue's dominated
    convergence theorem can be applied, so that we finally obtain that
$\theta\mapsto\var_\theta(\tilde\theta)$ is  c\`adl\`ag and fully continuous in those parameter values
that are not binary rationals.   

It remains to prove the formula for the jumps of $\phi$ in binary rationals.
Using the notation from the first part of the proof it is clear that $\theta\phi(\theta)$ is 
the finite sum $\sum_{m=1}^K a_m a_{m-1} 2^{-j_m}$, where $j_1,\dots,j_K=L$ are the successive 
positions of the digit~1 in the (finite) binary expansion of $\theta=q2^{-L}$. Similarly, 
$\theta\phi(\theta-)$ is the infinite sum $\sum_{m=1}^\infty a_m' a_{m-1}' 2^{-j_m'}$ based on
the sequence $(j_m')_{m\in\bN}$ with $j_m'=j_m$ for $m<K$, and $j'_{K-1+m}=L+m$ for all $m\in\bN$.
In particular, $a'_m=a_m$ for $m<K$, and we get
\begin{equation*}
	a'_{K-1+m} \,=\,  a'_{K-1}+2^{-L-1}+\cdots+2^{-L-m} \, = \, a_{K-1} + 2^{-L}(1-2^{-m})
	\  \text{ for all } m\in \bN_0,
\end{equation*}
which leads to
\begin{align*}
 \sum_{k=K}^\infty &a'_k a'_{k-1}(a'_k-a'_{k-1})\ 	
    	\ = \   \sum_{m=0}^\infty a'_{K+m} \,a'_{K+m-1} \, 2^{-j'_{K+m}}  \\
	      &= \   \sum_{m=0}^\infty \bigl(a_{K-1} + 2^{-L}(1-2^{-m-1}) \bigr)
	                    \bigl(a_{K-1} + 2^{-L}(1-2^{-m})\bigr)\, 2^{-L-m -1} \\
	      &=\ a_{K-1}^2 \, \sum_{m=0}^\infty  2^{-L - m -1}\\
	      &\hspace{1cm} + \ a_{K-1}\, \sum_{m=0}^\infty \bigl(2^{-L}(1-2^{-m-1}) 
	                              + \; 2^{-L}(1-2^{-m})  \bigr)\, 2^{-L - m -1}\\
	      &\hspace{1cm} +\   \sum_{m=0}^\infty     2^{-L}(1-2^{-m-1})  \,
                               	2^{-L}(1-2^{-m})\,2^{-L-m-1}         \\
	      &=\ a_{K-1}^2 \,2^{-L}   \,+\; a_{K-1}\, 2^{-2L} \;+\; \frac{2}{7}\,2^{-3L}, 
\end{align*} 
hence 
\begin{equation*}
	\theta \, \phi(\theta-)  
 	        = \, \sum_{k=1}^{K-1} a_k a_{k-1}(a_k-a_{k-1}) \; + \; a_{K-1}^2 \,2^{-L}\,+\; a_{K-1}\, 2^{-2L} 
                      	\;+\;\frac{2}{7}\,2^{-3L}.
\end{equation*}
From~\eqref{eq:MSEtilde} we know that
\begin{equation*}
	\theta \, \phi(\theta)  
 	        = \, \sum_{k=1}^{K-1} a_k a_{k-1}(a_k-a_{k-1}) \; + \;  a_K a_{K-1}(a_K-a_{K-1}).
\end{equation*}
Using $a_K=\theta$ and $a_K-a_{K-1}=2^{-L}$ we now obtain~\eqref{eq:sprung} after a short calculation.

\smallbreak
\textbf{Proof of Theorem~\ref{thm:aseff}.\ }
Fix $P_p\in\cP_0$. Denote by $\Theta_d$ the subset of elements $q=(q_k)_{k\in\bN_0}$ in $\Theta$
with finitely many non-zero mixing coefficients $q_k.$ 
For $q=(q_k)_{k\in\bN_0}$ in $\Theta_d$, $0\le a<\infty$ and $0\le t<\min(1,1/a)$ let
\begin{equation*}
   P_{p,q;a,t}=P_p+a t(P_q-P_p).
\end{equation*}
   We consider the submodel
\begin{equation*}
    \cP_{p,q;a}=\bigl\{P_{p,q;a,t}:\,  0\le t<\min(1,1/a)\bigr\}\; \subset\; \cP_0.
\end{equation*}
Then $\pp_p=\sum_{k=0}^\infty p_kf_{2k+1}$ is a density of $P_p$, and
\begin{equation*}
     \pp_{p,q;a,t}=\pp_p+at\sum_{k=0}^\infty (q_k-p_k)f_{2k+1}
\end{equation*}
is a density of $P_{p,q;a,t}$, both with respect to the Lebesgue measure $\Leb_+$ on
the Borel subsets of the positive half-line. 
The partial derivative
\[\dot{\pp}_{p,q;a,t}=\frac{\partial}{\partial t}\pp_{p,q;a,t}=a\sum_{k=0}^\infty (q_k-p_k)f_{2k+1}\]
does not depend on $t$, hence $t\rightarrow \pp_{p,\lambda;a,t}^{1/2}$ is continuously 
differentiable on the interval $[0,\min(1,1/a))$. Note that 
\[c_{p,q}:=\max \left (1,\frac{q_k}{p_k}, k\in\bN_0 \right )<\infty,\]
because only finitely many of the $q_k$ are non-zero. Due to  
\begin{align*}
\frac{\dot{\pp}_{p,q;a,t}^2}{\pp_{p,q;a,t}}
     \ &=\ a^2\frac{\left (\sum_{k=0}^\infty q_kf_{2k+1}-\sum_{k=0}^\infty p_kf_{2k+1}\right )^2}
          {(1-at)\sum_{k=0}^\infty p_kf_{2k+1}+at\sum_{k=0}^\infty q_kf_{2k+1}}\\
     \ &\le\ a^2\frac{\left (\sum_{k=0}^\infty\frac{q_k}{p_k}p_kf_{2k+1}+\sum_{k=0}^\infty p_kf_{2k+1}\right )^2}
         {(1-at)\sum_{k=0}^\infty p_kf_{2k+1}+at\sum_{k=0}^\infty q_kf_{2k+1}}\\
     \ &\le\ a^2\frac{(1+c_{p,q})^2}{1-at}\sum_{k=0}^\infty p_k f_{2k+1}(x),
\end{align*}                  
dominated convergence can be applied and it follows that the function 
$t\rightarrow \int \frac{\dot{\pp}_{p,q;a,t}^2}{\pp_{p,q;t}}\,d\Leb_+$
is continuous. Thus, $t\rightarrow P_{p,q;a,t}$ is a differentiable path with score function
\begin{equation*}
   g_{p,q;a}=\frac{\partial}{\partial t}\log \pp_{p,q;a,t}\Big |_{t=0}
      \; =\; 
   a \left[\frac{\sum_{k=0}^\infty q_kf_{2k+1}}{\sum_{k=0}^\infty p_k f_{2k+1}}-1 \right ],
\end{equation*}
meaning that
\begin{equation*}
  \bigintss \left [\frac{\pp_{\theta,q;a,t}^{1/2}-\pp_{p}^{1/2}}{t}-\frac{1}{2}g_{p,q;a}\hskip0.5mm \pp_p^{1/2}\right ]^2 d\Leb_+\rightarrow 0\ \text{ as }t\to 0,
\end{equation*}
see \cite[Lemma 1.8]{Vaart}. The tangent set
\begin{equation*}
  \dot{\cP}_{P_p}=\left \{ g_{p,q;a}:\, q\in \Theta_d,\, a\ge 0\right \}
\end{equation*}
of the model $\cP_0$ at $P_p$ is a convex cone. Regarding this cone as a subset of
$L_2(P_p)$, we claim that $\tilde{\kappa}_{P_p}$ defined by
\begin{equation*}
     \tilde{\kappa}_{P_p}(x)=\frac{1}{2}(x-1)-\kappa(P_p),\quad x >0,
\end{equation*}
is an element of the closed linear span $\overline{\rm lin}\dot{\cP}_{P_p}$ of $\dot{\cP}_{P_p}$.
To see this, we first note that with 
\begin{equation*}
  (2k-1)!!=\begin{cases}1&~\text{if}~k=0,\\
      (2k-1)\cdot (2k-3)\cdots 3\cdot 1 &~\text{if}~k\in \bN,
    \end{cases}
\end{equation*}
it holds that
\begin{equation*}
    f_{2k+1}(x)=\frac{1}{(2k-1)!!}x^k f_1(x),\hskip2mm k\in\bN_0.
\end{equation*}
Thus, for each $r\in\bN_0,$ and $\delta_r=(\delta_{rk})_{k\in\bN_0}$ with $\delta_{rr}=1$ and $\delta_{rk}=0$ otherwise,
\begin{equation*}
  g_{p,\delta_r;1}(x)\; =\; \frac{\frac{1}{(2r-1)!!}x^r}
                  {\sum_{k=0}^\infty \frac{1}{(2k-1)!!}p_kx^k}\, -\, 1,~x>0.
\end{equation*}
Putting
\begin{align*}
  a_0 & = -\left (\kappa(P_p)+\frac{1}{2}\right )p_0,\\
  a_k & = \frac{1}{2}(2k-1)p_{k-1} - \left (\kappa(P_p)+\frac{1}{2}\right )p_k,~k\in\bN,
\end{align*}
it follows that
\begin{equation*}
   \tilde{\kappa}_{P_p}(x) = \sum_{k=0}^\infty a_k g_{p,\delta_k;1}(x),~x>0,
\end{equation*}
where the series converges pointwise and in $L_2(P_p).$ Therefore, as asserted, $\tilde{\kappa}_{P_p}\in \overline{\rm lin}\dot{\cP}_{P_p}$. 
Due to
\begin{equation*}
      \kappa(P_{p,q;a,t})-\kappa(P_p)=t\int \tilde{\kappa}_{P_p}\, g_{p,q;a}\,dP_p  
\end{equation*}
for each $q\in \Theta_d$, each $a\ge 0$ and each $0\le t<\min(1,1/a)$, 
the functional $\kappa$ is differentiable at $P_p$ relative to the tangent cone $\dot{\cP}_{P_p}$.
Hence $\tilde{\kappa}_{P_p}$ is the efficient influence function for estimating the functional $\kappa$ at
$P_p$.  By \eqref{eq:unbiasedreg} and Lemma 2.9 in \cite{Vaart}, the sequence of estimators
$(T_n)_{n=1}^\infty$ for estimating $\kappa$ is thus asymptotically efficient at $P_p$.

\bibliographystyle{alpha}

\begin{thebibliography}{}
	
\bibitem{Amari} {\sc Amari, Shun-Ichi } (1982).
Differential geometry of curved exponential families---curvature and information loss. 
{\em Ann. Stat.} 10, 357-385.

\bibitem{BaGr5} {\sc Baringhaus, L. and Gr\"ubel, R.} (2020). 
  Mixture representations of noncentral distributions.  Published online in
  {\em Comm. Statist. Theory Methods}. DOI: 10.1080/03610926.2020.1738487.

\bibitem{Benedetto76} {\sc Benedetto, J.J.} (1976). {\em Real variable and integration.}
Teubner, Stuttgart.

\bibitem{DLR}{\sc Dempster, A. P., Laird, N. M. and Rubin, D. B.} (1977).
Maximum likelihood from incomplete data via the {EM} algorithm. 
{\em J. Roy. Statist. Soc. Ser. B} 39, 1--38.          

\bibitem{Dynkin78}{\sc Dynkin, E. B.} (1978). Sufficient statistics and extreme points. 
{\em Ann. Prob.} 6, 705--730.

\bibitem{Efron}{\sc Efron, B.} (1975). Defining the curvature of a statistical problem
  (with applications to second order efficiency). {\em Ann. Stat.} 3, 1189--1242.

\bibitem{FellerII}{\sc Feller, W.} (1971). {\em An Introduction to Probability Theory 
and Its Applications}, Vol. II, 2nd ed. Wiley, New York. 

\bibitem{Ferguson}{\sc Ferguson, T.S.} (1974).
Prior distributions on spaces of probability measures. 
{\em Ann. Statist.} 2, 615--629.

\bibitem{Ghosal-vdVaart}{\sc Ghosal, S. and van der Vaart, A.} (2017).
{\em Fundamentals of Nonparametric Bayesian Inference.}
Cambridge University Press, Cambridge.

\bibitem{GoelDeGroot79}{\sc Goel, P.K. and DeGroot, M.H.} (1979). Comparison of experiments and information measures. {\em Ann. Statist.} 7, 1066--1077. 

\bibitem{Grenander}  {\sc Grenander, Ulf} (1956). 
On the theory of mortality measurement II. {\em Skand. Aktuarietidskr.} 39, 125--153.

\bibitem{Hoff} {\sc Hoff, Peter D.} (2003).
  Nonparametric estimation of convex models via mixtures,  {\em Ann. Statist.} 31, 174--200.

\bibitem{Kubokawa17}{\sc Kubokawa, T.,  Marchand, \'E. and Strawderman, W. E.} (2017).
A unified approach to estimation of noncentrality parameters, the multiple correlation coefficient,
and mixture models. {\em Math. Methods Statist.} 26, 134--148.           

\bibitem{Lauritzen} {\sc Lauritzen, S.} (1988). Extremal families and systems of sufficient statistics.
{\em Lecture Notes in Statistics} 49. Springer, Heidelberg.

\bibitem{LeCam96}{\sc Le Cam, L.} (1996). Comparison of experiments - a short review.
{\em Statistics, probability and game theory}, 127--138. IMS Lecture Notes Monogr. Ser., 30.
Inst. Math. Statist. Hayward, CA.  

\bibitem{Lehmann}{\sc Lehmann, E. L. and Casella, G.} (1998). {\em Theory of Point Estimation}, second edition.
  Springer-Verlag, New York. 

\bibitem{Leroux}{\sc Leroux, Brian G.} (1992).
Consistent estimation of a mixing distribution, {\em Ann. Statist.} 20, 1350--1360.

\bibitem{Lindsay} {\sc Lindsay, B. G.} (1995). {\em Mixture Models: Theory, Geometry 
and Applications.} NSF-CBMS Regional Conference Series in Probability and Statistics, 
Volume 5. IMS, Hayward.

\bibitem{pfanzagl} {\sc Pfanzagl, J.} (1988). Consistency of maximum likelihood estimators
  for certain nonparametric families, in particular: mixtures. {\em J. Statist. Plann. Inference} 19, 137--158.

\bibitem{Phelps}{\sc Phelps, Robert R.} (2001). {\em Lectures on {C}hoquet's theorem}, second edition. Lecture Notes in Mathematics 1757, Springer-Verlag, Berlin.

\bibitem{RudinFA}{\sc Rudin, W.} (1991). {\em Functional Analysis.} McGraw-Hill, New York.

\bibitem{SaxenaAlam82}{\sc Saxena, K. M. L. and Alam, K.} (1982). Estimation of the noncentrality parameter of a chi squared distribution. {\em Ann. Statist.} 10, 1012--1016.           
  
\bibitem{Stone61}{\sc Stone, M.} (1961). Non-equivalent comparisons of experiments and
  their use for experiments involving location parameters. {\em Ann. Math. Statist.} 32, 326--332.

\bibitem{Vaart} {\sc van der Vaart, A.} (2002). Semiparametric statistics. 
In: Bolthausen, E., Perkins, E. and van der Vaart, A. {\em Lectures on Probability 
Theory and Statistics.} Lecture Notes in Mathematics 1781, Springer, Berlin.

\bibitem{WoodProb} {\sc Wood, G. R.} (1992).
Binomial mixtures and finite exchangeability, {\em Ann. Probab.} 20, 1167--1173.

\bibitem{WoodStat}{\sc Wood, G. R.} (1999).
Binomial mixtures: geometric estimation of the mixing distribution,
{\em Ann. Statist.} 27, 1706--1721.
	
\end{thebibliography}


\end{document}